\pdfoutput=1
\RequirePackage{ifpdf}
\ifpdf % We~are running pdfTeX in pdf mode
\documentclass[pdftex]{sigma}
\else
\documentclass{sigma}
\fi

\numberwithin{equation}{section}

\newtheorem{Theorem}{Theorem}[section]
\newtheorem*{Theorem*}{Theorem}
\newtheorem{Corollary}[Theorem]{Corollary}
\newtheorem{Lemma}[Theorem]{Lemma}
\newtheorem{Proposition}[Theorem]{Proposition}

\newtheorem{mainthm}{Theorem}

\theoremstyle{definition}
\newtheorem{Definition}[Theorem]{Definition}

\usepackage{dsfont}
\usepackage{stmaryrd}

\usepackage[scr=boondox]{mathalfa}

\newcommand{\Ad}{\mathrm{Ad}}
\DeclareMathOperator{\Hol}{\mathrm{Hol}}
\DeclareMathOperator{\Aut}{\mathrm{Aut}}
\DeclareMathOperator{\Iso}{\mathrm{Euc}}
\DeclareMathOperator{\Orth}{\mathrm{O}}
\DeclareMathOperator{\Utry}{\mathrm{U}}

\DeclareMathOperator{\GLin}{\mathrm{GL}}

\DeclareMathOperator{\PGL}{\mathrm{PGL}}
\DeclareMathOperator{\PU}{\mathrm{PU}}

\DeclareMathOperator{\Path}{Path}

\newcommand{\MC}[1]{\omega_{{}_{#1}}}
\newcommand{\Lt}[1]{\operatorname{L}_{#1}}
\newcommand{\Rt}[1]{\operatorname{R}_{#1}}

\usepackage[all]{xy}

\usepackage[normalem]{ulem}

\begin{document}

\allowdisplaybreaks

\newcommand{\arXivNumber}{2406.05497}

\renewcommand{\PaperNumber}{084}

\FirstPageHeading

\ShortArticleName{Holonomy of Parabolic Geometries Near Isolated Higher-Order Fixed Points}

\ArticleName{Holonomy of Parabolic Geometries Near Isolated \\ Higher-Order Fixed Points}

\Author{Jacob W. ERICKSON}

\AuthorNameForHeading{J.W.~Erickson}

\Address{School of Mathematics, Georgia Institute of Technology, Atlanta, USA}
\Email{\mail{jerickson41@gatech.edu}}
\URLaddress{\url{https://erickson.math.gatech.edu/}}

\ArticleDates{Received December 01, 2025, in final form August 12, 2026; Published online August 28, 2026}

\Abstract{For Cartan geometries admitting an automorphism with isotropy satisfying a particular, loosely dynamical property on their model geometries, we demonstrate the existence of an open subset of the geometry with trivial holonomy. This property, which generalizes characteristics of isotropies corresponding to isolated higher-order fixed points in parabolic geometries that are known to require a nearby open subset to have vanishing curvature, only relies upon the behavior of the isotropy in the model geometry, and therefore applies regardless of initial curvature assumptions, such as regularity or normality. Along the way to proving our main results, we also derive a couple of results for working with holonomy, relating to limits of sequences of developments and the existence of antidevelopments, that are useful in their own right. To showcase the effectiveness of the techniques developed, we use them to completely characterize all almost c-projective and almost quaternionic structures that admit a nontrivial automorphism with a higher-order fixed point, as well as all nondegenerate partially integrable almost CR structures that admit a higher-order fixed point with non-null isotropy.}


\Classification{53C29; 58D19; 57S25; 53C30}

\section{Introduction}
In each dimension greater than two, a celebrated theorem of Ferrand and Obata tells us that, up to conformal isomorphism, there are only two examples of conformal structures of Riemannian signature admitting automorphism groups that act non-properly: the standard conformal sphere and Euclidean space. A similar result, due to Schoen and Webster, holds for strictly pseudo-convex CR structures: if the CR automorphisms of such a structure act non-properly, then it must be isomorphic to either the standard CR structure on the sphere or a left-invariant CR structure on the Heisenberg group. Eventually, it was realized that Cartan geometries---a class of geometric structures that extend the spirit of Klein's Erlangen program somewhat analogously to the way that Riemannian manifolds extend Euclidean geometry---are the natural setting for this type of result, and in \cite{Frances2007}, Frances gave a generalization of both results, extending the idea to all parabolic Cartan geometries of real rank one, assuming the curvature satisfies a standard regularity condition. A convenient and concise history of this result can be found in \cite{MelnickFO}.

\!These Ferrand--Obata type theorems provide striking examples where the global structure of a given Cartan geometry is almost completely determined by the behavior of its automorphisms. Naturally, though, we are led to ask: how far can we take this? How much can we say about the global structure of a Cartan geometry just from the behavior of its automorphisms, or even from the behavior of a \textit{single} automorphism?

In the vein of this latter question, we specifically want to look at automorphisms of parabolic geometries with \textit{higher-order fixed points}. These are automorphisms with isotropies, over some fixed point of the geometry, lying in the unipotent radical of the parabolic subgroup used as the stabilizer subgroup for the corresponding model geometry. For many of the most common types of parabolic geometries, such as projective and conformal geometry, these just amount to automorphisms whose derivatives on the base manifold are the identity transformation on the tangent space over some fixed point.

The existence of a higher-order fixed point often guarantees that the curvature of the geometry vanishes over some open subset containing the higher-order fixed point in its closure. Results in this direction essentially go all the way back to \cite{NaganoOchiai} in the projective case and \cite{Frances2012} in the conformal case, though a more unified approach was introduced in \cite{CapMelnick2013} and expanded upon in~\cite{KruglikovThe2018,MelnickNeusser2016}, among other places. Notably, however, these prior results are all local in nature, unless we make quite strong assumptions on the geometry, such as real-analyticity or---as explored in~\cite{TianyuMa, NaganoOchiai}---metrizability.

In this paper, we will provide a way to extend these existing local results to semi-global ones when the higher-order fixed point is isolated. Our main result, Theorem~\ref{ensnaringthm}, proves that if an automorphism has isotropy satisfying a kind of dynamical property in the model geometry that applies to the cases with isolated higher-order fixed points where local results are known, then there is a certain invariant open subset of the geometry on which the \emph{holonomy} must be trivial. If we know a bit more about the topology of this open subset, then we can even prove genuine global results along the lines of Ferrand--Obata type theorems. To showcase this, we specifically prove such results for c-projective structures, quaternionic structures, and CR structures; as a~courtesy to those unfamiliar with Cartan geometries, we state these theorems in their ``native'' languages, though in the rest of the paper, we will almost exclusively work from the unified perspective provided by the Cartan machinery.

\begin{mainthm}\label{thmA} Suppose that $(M,J,[\nabla])$ is an almost c-projective structure.
If $(M,J,[\nabla])$ admits a nontrivial automorphism $\alpha$ with a higher-order fixed point $x\in M$, so that $\alpha(x)=x$ and $\alpha_{*x}=\mathrm{id}_{T_xM}$, then $(M,J,[\nabla])$ embeds geometrically onto a dense open subspace of the standard flat c-projective structure on $\mathbb{CP}^m$. In particular, if $M$ is compact, then such an automorphism and fixed point exist if and only if $(M,J,[\nabla])$ is isomorphic to the flat c-projective structure on~$\mathbb{CP}^m$.\end{mainthm}

\begin{mainthm}\label{thmB} Suppose that $M$ is an almost quaternionic manifold. If $M$ admits a nontrivial automorphism $\alpha$ with a higher-order fixed point $x\in M$, so that $\alpha(x)=x$ and $\alpha_{*x}=\mathrm{id}_{T_xM}$, then $M$ embeds geometrically onto a dense open subspace of $\mathbb{HP}^m$ equipped with its standard quaternionic structure. In particular, if $M$ is compact, then such an automorphism and fixed point exist if and only if $M$ is isomorphic to $\mathbb{HP}^m$.\end{mainthm}

\begin{mainthm}\label{thmC} Suppose that $(M,J,D)$ is a partially integrable almost CR structure with Levi form of signature $(p,q)$, where $J$ is an almost complex structure on the contact distribution $D$.
If $(M,J,D)$ admits a nontrivial automorphism $\alpha$ with a non-null higher-order fixed point $x\in M$, so that $\alpha(x)=x$ and either $\alpha_{*x}=\mathrm{id}_{T_xM}$ or the image of the difference $\alpha_{*x}-\mathrm{id}_{T_xM}$ consists of non-null vectors for the Levi form on the contact distribution $D$, then $(M,J,D)$ embeds geometrically onto a dense open subspace of $\mathrm{Null}(\mathrm{h}_{p,q}):=\bigl\{\mathbb{C}^\times u\in\mathbb{CP}^{p+q+1}\mid \mathrm{h}_{p,q}(u,u)=0\bigr\}$ equipped with the standard flat CR structure, where the Hermitian form $\mathrm{h}_{p,q}$ on $\mathbb{C}^{p+q+2}$ has corresponding quadratic form given by \[\begin{bmatrix}z_0 \\ \vdots \\ z_{p+q+1}\end{bmatrix}\mapsto 2\mathrm{Re}(\bar{z}_0z_{p+q+1})+\sum_{j=1}^p|z_j|^2-\sum_{j=p+1}^{p+q}|z_j|^2.\] In particular, if $M$ is compact, then such an automorphism and fixed point exist if and only if $(M,J,D)$ is isomorphic to the flat CR structure on $\mathrm{Null}(\mathrm{h}_{p,q})$.\end{mainthm}

These results will follow from Theorem~\ref{ensnaringthm+}, which is itself a corollary of our main result, Theorem~\ref{ensnaringthm}. To get to these results, however, we will need to prove a couple of intermediate facts that are useful for working with development and holonomy in this setting, which we provide in Section~\ref{devtools} after reviewing terminology and notation in Section~\ref{preliminaries}. Then, in Section~\ref{snares}, we will describe the dynamical property for isotropies alluded to above and use it to prove our main result, before finishing the paper with the applications in Section~\ref{applications}.

We should remark that Theorems \ref{thmA}, \ref{thmB}, and \ref{thmC} all appeared in the original preprint for~\cite{SprawlPaper}. This current paper and its methods grew, at least in part, out of efforts to simplify the proofs and assumptions in that preprint, to avoid having to use the somewhat heavy-handed sprawl construction. It should be considered as separate from---though obviously related to---the current version of~\cite{SprawlPaper}, which is substantially different and does not include these results.

\section{Preliminaries}\label{preliminaries}
This section is meant as a way of establishing terminology and notation; it is probably not a good introduction to the topic. We recommend \cite{Sharpe1997} for an actual introduction to Cartan geometries, \cite{CapSlovakPG1} for an introduction to parabolic geometries, and \cite{HolonomyPaper} for an overview of the results on holonomy.

\subsection{Relevant model geometries}
Cartan geometries are modeled on homogeneous geometries in the sense of Klein, though we present these geometries in a way that emphasizes the role of the Lie group as a principal bundle over the homogeneous space.

\begin{Definition}A \emph{model} (or \emph{model geometry}) is a pair $(G,H)$, where $G$ is a Lie group and $H$ is a closed subgroup of $G$ such that $G/H$ is connected. In this case, $G$ is called the \emph{model group} and $H$ is called the \emph{isotropy} or \emph{stabilizer subgroup}.\end{Definition}

Writing $\Iso(m):=\mathbb{R}^m\rtimes\Orth(m)$ for the group of Euclidean isometries on $\mathbb{R}^m$ (with the dot product metric structure) and $\mathrm{Aff}(m):=\mathbb{R}^m\rtimes\GLin_m\mathbb{R}$ for the group of affine transformations on $\mathbb{R}^m$, two standard examples of model geometries are $(\Iso(m),\Orth(m))$, which corresponds to Euclidean geometry on $\mathbb{R}^m\cong\Iso(m)/\Orth(m)$, and $(\mathrm{Aff}(m),\GLin_m\mathbb{R})$, which corresponds to affine geometry on $\mathbb{R}^m\cong\mathrm{Aff}(m)/\GLin_m\mathbb{R}$. Cartan geometries that are modeled on these two examples correspond to Riemannian structures (modulo torsion considerations) and affine connections, respectively.

The models relevant to Theorems \ref{thmA}, \ref{thmB}, and \ref{thmC} are all \emph{parabolic}.

\begin{Definition}A model $(G,P)$ is \emph{parabolic} whenever the model group $G$ is semisimple and the isotropy\footnote{When we assume that the model is parabolic, we will always denote the isotropy subgroup by $P$. This will help us to visually distinguish between results that apply to more general Cartan geometries and results that apply to particular parabolic geometries.} $P$ is a parabolic subgroup.\end{Definition}

For parabolic models, we get an $\Ad_P$-invariant filtration of $\mathfrak{g}$ given by \[\mathfrak{g}=\mathfrak{g}^{-k}\supset\mathfrak{g}^{-k+1}\supset\cdots\supset\mathfrak{g}^k\supset\{0\},\] with $\mathfrak{g}^0:=\mathfrak{p}$ and $\mathfrak{g}^1:=\mathfrak{p}_+$, where $\mathfrak{p}_+$ is the nilradical of $\mathfrak{p}$. From this filtration, we can also get a grading $\sum_{i=-k}^k\mathfrak{g}_i$ on $\mathfrak{g}$ with $\mathfrak{g}_i\approx\mathfrak{g}^i/\mathfrak{g}^{i+1}$; this grading corresponds to a choice of Cartan involution, but we will keep this choice implicit, since we will not need the Cartan involution here. Notably, the filtration and grading satisfy $[\mathfrak{g}^i,\mathfrak{g}^j]\subseteq\mathfrak{g}^{i+j}$ and $[\mathfrak{g}_i,\mathfrak{g}_j]\subseteq\mathfrak{g}_{i+j}$ for all~$i$ and~$j$, so $\mathfrak{g}_-:=\sum_{i<0}\mathfrak{g}_i$ and $\mathfrak{g}_0$ are subalgebras. We denote by $G_-$ and $P_+$ the connected nilpotent subgroups generated by $\mathfrak{g}_-$ and $\mathfrak{p}_+$, respectively, and by $G_0$ the closed subgroup of $P$ such that $\Ad_{G_0}(\mathfrak{g}_i)=\mathfrak{g}_i$ for every grading component~$\mathfrak{g}_i$.

In Theorems \ref{thmA}, we will focus on the model $(\PGL_{m+1}\mathbb{C},P)$, where $\PGL_{m+1}\mathbb{C}$ is the quotient of $\GLin_{m+1}\mathbb{C}$ by its center $\mathbb{C}^\times\mathds{1}$ and \[P:=\left\{\begin{pmatrix}r & p \\ 0 & A\end{pmatrix}\in\PGL_{m+1}\mathbb{C}\mid r\in\mathbb{C}^\times,\,p^\top\in\mathbb{C}^m,\, A\in\GLin_m\mathbb{C}\right\}.\] This model geometry is parabolic, with grading
\begin{gather*}
\mathfrak{g}_-=\mathfrak{g}_{-1}:=\left\{\begin{pmatrix}0 & 0 \\ v & 0\end{pmatrix}\in\mathfrak{pgl}_{m+1}\mathbb{C}\mid v\in\mathbb{C}^m\right\},\\
\mathfrak{g}_0:=\left\{\begin{pmatrix}r & 0 \\ 0 & R\end{pmatrix}\in\mathfrak{pgl}_{m+1}\mathbb{C}\mid r\in\mathbb{C},\, R\in\mathfrak{gl}_m\mathbb{C}\right\},
\end{gather*}
 and \[\mathfrak{p}_+=\mathfrak{g}_1:=\left\{\begin{pmatrix}0 & p \\ 0 & 0\end{pmatrix}\in\mathfrak{pgl}_{m+1}\mathbb{C}\mid p^\top\in\mathbb{C}^m\right\},\] and associated subgroups
 \begin{gather*}
 G_-:=\left\{\begin{pmatrix}1 & 0 \\ v & \mathds{1}\end{pmatrix}\in\PGL_{m+1}\mathbb{C}\mid v\in\mathbb{C}^m\right\},\\ G_0:=\left\{\begin{pmatrix}r & 0 \\ 0 & A\end{pmatrix}\in\PGL_{m+1}\mathbb{C}\mid r\in\mathbb{C}^\times,\,  A\in\GLin_m\mathbb{C}\right\},
  \end{gather*}
  and \[P_+:=\left\{\begin{pmatrix}1 & p \\ 0 & \mathds{1}\end{pmatrix}\in\PGL_{m+1}\mathbb{C}\mid p^\top\in\mathbb{C}^m\right\}.\] In this case, the model encodes complex projective geometry over $\PGL_{m+1}\mathbb{C}/P\cong\mathbb{CP}^m$, and Cartan geometries of type $(\PGL_{m+1}\mathbb{C},P)$ (with certain curvature restrictions) correspond to almost c-projective structures. An overview of such structures can be found in~\cite{cprojCEMN}.

In Theorem~\ref{thmB}, we will focus on the model $(\PGL_{m+1}\mathbb{H},P)$. Here, $\PGL_{m+1}\mathbb{H}$ is the quotient of the quaternionic general linear group $\GLin_{m+1}\mathbb{H}$ of right $\mathbb{H}$-module automorphisms of $\mathbb{H}^{m+1}$ by its center, which corresponds to those automorphisms which left-multiply every element of~$\mathbb{H}^{m+1}$ by a nonzero real number, and \[P:=\left\{\begin{pmatrix}r & p \\ 0 & A\end{pmatrix}\in\PGL_{m+1}\mathbb{H}\mid r\in\mathbb{H}^\times,\, p^\top\in\mathbb{H}^m,\, A\in\GLin_m\mathbb{H}\right\}.\] This model geometry is also parabolic, with grading, subgroups, and subalgebras the same as those for $(\PGL_{m+1}\mathbb{C},P)$ but with~$\mathbb{C}$ replaced by $\mathbb{H}$ everywhere. It encodes the quaternionic analogue of projective geometry over $\PGL_{m+1}\mathbb{H}/P\cong\mathbb{HP}^m$, and Cartan geometries of type $(\PGL_{m+1}\mathbb{H},P)$ (again with certain curvature restrictions) correspond to almost quaternionic structures, as described in \cite[Proposition~4.1.8]{CapSlovakPG1}.

Finally, in Theorem~\ref{thmC}, we will focus on the model $(\PU(\mathrm{h}_{p,q}),P)$, where $\mathrm{h}_{p,q}$ is the Hermitian form on $\mathbb{C}^{p+q+2}$ with quadratic form given by \[\begin{bmatrix}z_0 \\ \vdots \\ z_{p+q+1}\end{bmatrix}\mapsto 2\mathrm{Re}(\bar{z}_0z_{p+q+1})+\sum_{j=1}^p |z_j|^2-\sum_{j=p+1}^{p+q}|z_j|^2,\] and $\PU(\mathrm{h}_{p,q})$ is the quotient of the group of unitary transformations for $\mathrm{h}_{p,q}$ by its center, consisting of multiples of the identity matrix by elements of $\Utry(1)$. Denoting by $I_{p,q}$ the ${(p+q)}\times(p+q)$ diagonal matrix with the first $p$ diagonal entries equal to $+1$ and the last $q$ diagonal entries equal to $-1$, the Lie group $\PU(\mathrm{h}_{p,q})$ has Lie algebra of the form
\[\mathfrak{pu}(\mathrm{h}_{p,q}):=\left\{\begin{pmatrix}r & \beta & \mathrm{i}s \\ v & R & -I_{p,q}\bar{\beta}^\top \\ \mathrm{i}t & -\bar{v}^\top I_{p,q} & -\bar{r}\end{pmatrix}\,\Big| \begin{array}{l@{}} s,t\in\mathbb{R},\, v,\beta^\top\in\mathbb{C}^{p+q}, \\ r\in\mathbb{C},\text{ and } R\in\mathfrak{u}(p,q)\end{array}\right\},\]
where elements of the Lie algebra are considered equivalent if their difference is an imaginary multiple of the identity matrix. The parabolic subgroup $P$, then, is
\[P:=\left\{\begin{pmatrix}r & r\beta & r\bigl(\mathrm{i}s-\tfrac{1}{2}\beta I_{p,q}\bar{\beta}^\top\bigr) \vspace{1mm}\\ 0 & A & -AI_{p,q}\bar{\beta}^\top \vspace{1mm}\\ 0 & 0 & \bar{r}^{-1}\end{pmatrix}\in\PU(\mathrm{h}_{p,q})\,\Big| \begin{array}{l@{}} s\in\mathbb{R},\, \beta^\top\in\mathbb{C}^{p+q}, \\ r\in\mathbb{C}^\times,\text{ and }A\in\Utry(p,q)\end{array}\right\},\]
with grading given by
\begin{gather*}
\mathfrak{g}_{-2}:=\left\langle\begin{pmatrix}0 & 0 & 0 \\ 0 & 0 & 0 \\ \mathrm{i} & 0 & 0\end{pmatrix}\right\rangle,\\
 \mathfrak{g}_{-1}:=\left\{\begin{pmatrix}0 & 0 & 0 \\ v & 0 & 0 \\ 0 & -\bar{v}^\top I_{p,q} & 0\end{pmatrix}\in\mathfrak{pu}(\mathrm{h}_{p,q})\mid v\in\mathbb{C}^{p+q}\right\},\\
 \mathfrak{g}_0:=\left\{\begin{pmatrix}r & 0 & 0 \\ 0 & R & 0 \\ 0 & 0 & -\bar{r}\end{pmatrix}\in\mathfrak{pu}(\mathrm{h}_{p,q})\mid r\in\mathbb{C},\, R\in\mathfrak{u}(p,q)\right\},\\
 \mathfrak{g}_1:=\left\{\begin{pmatrix}0 & \beta & 0 \\ 0 & 0 & -I_{p,q}\bar{\beta}^\top \\ 0 & 0 & 0\end{pmatrix}\in\mathfrak{pu}(\mathrm{h}_{p,q}) \mid \beta^\top\in\mathbb{C}^{p+q}\right\},
\end{gather*}
and
\[\mathfrak{g}_2:=\left\langle\begin{pmatrix}0 & 0 & \mathrm{i} \\ 0 & 0 & 0 \\ 0 & 0 & 0\end{pmatrix}\right\rangle,\]
with corresponding subgroups
\begin{gather*}
G_-:=\left\{\begin{pmatrix}1 & 0 & 0 \\ v & \mathds{1} & 0 \\ \mathrm{i}t-\tfrac{1}{2}\bar{v}^\top I_{p,q} v & -\bar{v}^\top I_{p,q} & 1\end{pmatrix}\in\PU(\mathrm{h}_{p,q})\mid v\in\mathbb{C}^{p+q},\, t\in\mathbb{R}\right\},\\
 G_0:=\left\{\begin{pmatrix}r & 0 & 0 \\ 0 & A & 0 \\ 0 & 0 & \bar{r}^{-1}\end{pmatrix}\in\PU(\mathrm{h}_{p,q})\mid r\in\mathbb{C}^\times,\, A\in\Utry(p,q)\right\},
 \end{gather*}
and
\[P_+:=\left\{\begin{pmatrix}1 & \beta & \mathrm{i}s-\tfrac{1}{2}\beta I_{p,q} \bar{\beta}^\top \vspace{1mm}\\ 0 & \mathds{1} & -I_{p,q}\bar{\beta}^\top \\ 0 & 0 & 1\end{pmatrix}\in\PU(\mathrm{h}_{p,q})\mid \beta^\top\in\mathbb{C}^{p+q},\, s\in\mathbb{R}\right\}.\]
Within $P_+$, we will also say that \[a=\begin{pmatrix}1 & \beta & \mathrm{i}s-\tfrac{1}{2}\beta I_{p,q} \bar{\beta}^\top \vspace{1mm}\\ 0 & \mathds{1} & -I_{p,q}\bar{\beta}^\top \\ 0 & 0 & 1\end{pmatrix}\in P_+\] is \emph{non-null} if and only if it is nontrivial and $\beta I_{p,q}\bar{\beta}^\top\neq 0$ if $\beta\neq 0$.

The homogeneous space $\PU(\mathrm{h}_{p,q})/P$ is naturally diffeomorphic to the null-cone \[\mathrm{Null}(\mathrm{h}_{p,q}):=\bigl\{\mathbb{C}^\times u\in\mathbb{CP}^{p+q+1}\mid \mathrm{h}_{p,q}(u,u)=0\bigr\}\] for $\mathrm{h}_{p,q}$ in $\mathbb{CP}^{p+q+1}$. This null-cone is a compact and simply connected smooth manifold of (real) dimension $2(p+q)+1$; for $q=0$ or $p=0$, $\mathrm{Null}(\mathrm{h}_{p,q})$ is diffeomorphic to the sphere. Cartan geometries of type $(\PU(\mathrm{h}_{p,q}),P)$ (again, with certain curvature restrictions) correspond to nondegenerate partially integrable almost CR structures with Levi form of signature $(p,q)$, as detailed in \cite[Section~4.2.4]{CapSlovakPG1}.

\subsection{Generalities for Cartan geometries}
In essence, the idea of a Cartan geometry of type $(G,H)$ is to specify a $\mathfrak{g}$-valued one-form $\omega$ that behaves like the Maurer-Cartan form $\MC{G}\colon X_g\in T_g G\mapsto\Lt{g^{-1}*}X_g\in\mathfrak{g}$ does on $G$, where $\Lt{a}\colon g\mapsto ag$ denotes left-translation by $a\in G$.

\begin{Definition}Let $(G,H)$ be a model. A \emph{Cartan geometry of type $(G,H)$ over a (smooth) manifold $M$} is a pair $(\mathscr{G},\omega)$, where $\mathscr{G}$ is a principal $H$-bundle over $M$ with quotient map\footnote{In an effort to declutter notation, we will always denote the quotient map of a principal $H$-bundle by $q_{{}_H}$, even if there are multiple relevant principal bundles. The meaning of the map should always be clear from context. Similarly, we will always denote the right-action of $h\in H$ on a principal bundle by $\Rt{h}\colon \mathscr{g}\mapsto\mathscr{g}h$.} $q_{{}_H}\colon \mathscr{G}\to M$ and $\omega$ is a $\mathfrak{g}$-valued one-form on $\mathscr{G}$ such that
\begin{itemize}\itemsep=0pt
\item for every $\mathscr{g}\in\mathscr{G}$, $\omega_\mathscr{g} \colon T_\mathscr{g}\mathscr{G}\to\mathfrak{g}$ is a linear isomorphism;
\item for every $h\in H$, $\Rt{h}^*\omega=\Ad_{h^{-1}}\omega$;
\item for every $Y\in\mathfrak{h}$, the flow of the vector field $\omega^{-1}(Y)$ is given by $\exp\bigl(t\omega^{-1}(Y)\bigr)=\Rt{\exp(tY)}$ for all $t\in\mathbb{R}$.\end{itemize}\end{Definition}

A natural example of a Cartan geometry of type $(G,H)$ is always the \emph{Klein geometry} of that type, which encodes the geometric structure of the model geometry as a Cartan geometry.

\begin{Definition}For a model $(G,H)$, the \emph{Klein geometry of type $(G,H)$} is the Cartan geometry $(G,\MC{G})$ over $G/H$, where $G$ is the model group and $\MC{G}$ is the Maurer-Cartan form on $G$.\end{Definition}

Throughout, we will want to compare different Cartan geometries of the same type. To do this, we will use local diffeomorphisms that preserve the Cartan-geometric structure, called \emph{geometric maps}.

\begin{Definition}Given two Cartan geometries $(\mathscr{G},\omega)$ and $(\mathscr{Q},\upsilon)$ of type $(G,H)$, a \emph{geometric map} $\varphi\colon (\mathscr{G},\omega)\to (\mathscr{Q},\upsilon)$ is an $H$-equivariant smooth map $\varphi\colon \mathscr{G}\to\mathscr{Q}$ such that $\varphi^*\upsilon=\omega$.\end{Definition}

A geometric map $\varphi\colon (\mathscr{G},\omega)\to (\mathscr{Q},\upsilon)$ always induces a corresponding local diffeomorphism between the base manifolds of $\mathscr{G}$ and $\mathscr{Q}$, given by $q_{{}_H}(\mathscr{g})\mapsto q_{{}_H}(\varphi(\mathscr{g}))$ for each $\mathscr{g}\in\mathscr{G}$. We find it convenient and natural to not waste symbols to distinguish between these maps; throughout, whenever we have a geometric map $\varphi$, we will also denote its induced map on the base manifolds by the same symbol $\varphi$. The meaning should always be clear from context.

Of course, some geometric maps tell us more than others. We say that a geometric map $\varphi\colon (\mathscr{G},\omega)\to (\mathscr{Q},\upsilon)$ is a \emph{geometric embedding} when $\varphi$ is injective, and when $\varphi$ is bijective, we further say that it is a \emph{(geometric) isomorphism}. A geometric isomorphism from $(\mathscr{G},\omega)$ to itself is then called a \emph{(geometric) automorphism}.

Geometric maps are quite rigid: if the base manifold of $(\mathscr{G},\omega)$ is connected, then a geometric map $\varphi\colon (\mathscr{G},\omega)\to(\mathscr{Q},\upsilon)$ is uniquely determined by the image $\varphi(\mathscr{g})$ of an arbitrary single element $\mathscr{g}\in\mathscr{G}$. As such, the \emph{automorphism group} $\Aut(\mathscr{G},\omega)$, consisting of all automorphisms of $(\mathscr{G},\omega)$, will act freely on $\mathscr{G}$.

Following our convention of not distinguishing between geometric maps and the corresponding induced maps on the base manifolds, when we talk about fixed points of an automorphism $\alpha\in\Aut(\mathscr{G},\omega)$, we will mean fixed points of the induced map on the base manifold; since $\Aut(\mathscr{G},\omega)$ acts freely, automorphisms will not have actual fixed points in the bundle $\mathscr{G}$. We~denote by \[\mathrm{Fix}_M(\alpha):=\{q_{{}_H}(\mathscr{g})\in M \mid \alpha(q_{{}_H}(\mathscr{g}))=q_{{}_H}(\mathscr{g})\}\] the set of fixed points for $\alpha\in\Aut(\mathscr{G},\omega)$.

For each $\mathscr{g}\in\mathscr{G}$ such that $q_{{}_H}(\mathscr{g})\in\mathrm{Fix}_M(\alpha)$, there is some $a\in H$ such that $\alpha(\mathscr{g})=\mathscr{g}a$. This element $a\in H$ is called the \emph{isotropy} of $\alpha$ at $\mathscr{g}$.

\begin{Definition}For $\alpha\in\Aut(\mathscr{G},\omega)$ and $\mathscr{e}\in\mathscr{G}$ such that $\alpha(\mathscr{e})\in\mathscr{e}H$, the \emph{isotropy} of $\alpha$ at $\mathscr{e}$ is the unique element $a\in H$ such that $\alpha(\mathscr{e})=\mathscr{e}a$.\end{Definition}

For the applications considered in the paper, we will primarily focus on automorphisms $\alpha\in\Aut(\mathscr{G},\omega)$ of parabolic geometries $(\mathscr{G},\omega)$ of type $(G,P)$ with isotropy $a\in P_+<P$ at some $\mathscr{e}\in\mathscr{G}$. In that case, we say that $\alpha$ has a \emph{higher-order fixed point} at $q_{{}_P}(\mathscr{e})$.

Another core idea for Cartan geometries is that of \emph{curvature}, which tells us when the geometry locally differs from the Klein geometry.

\begin{Definition}Given a Cartan geometry $(\mathscr{G},\omega)$ of type $(G,H)$, its \emph{curvature} is the $\mathfrak{g}$-valued two-form given by $\Omega:=\mathrm{d}\omega+\tfrac{1}{2}[\omega,\omega]$.\end{Definition}

When the curvature vanishes in a neighborhood of a point, then the geometry is locally equivalent to that of the Klein geometry near that point. In other words, when $\Omega$ vanishes on some neighborhood of an element $\mathscr{e}\in\mathscr{G}$, we can find a geometric embedding \[\psi\colon \ \bigl(q_{{}_H}^{-1}(U),\MC{G}\bigr)\hookrightarrow (\mathscr{G},\omega)\] from a neighborhood $q_{{}_H}^{-1}(U)$ of the identity $e\in G$ in the Klein geometry such that $\psi(e)=\mathscr{e}$. When the curvature vanishes everywhere, we say that the geometry is \emph{flat}; such geometries are equivalent to locally homogeneous geometric structures in the sense of~\cite{Goldman2022}.

\subsection{Development and holonomy}
Again, we recommend \cite{HolonomyPaper} for an overview of our techniques involving holonomy, as well as \cite[Chapter~3, Section~7]{Sharpe1997} for a review of basic results on developments of paths.

We would like to pretend that Cartan geometries $(\mathscr{G},\omega)$ of type $(G,H)$ ``are'' their model geometries. The notions of \emph{development} and \emph{holonomy} allow us to do this somewhat judiciously.

\begin{Definition}Given a (piecewise smooth)\footnote{Throughout, whenever we refer to a ``path'', we will always mean a piecewise smooth path.} path $\gamma\colon [0,1]\to\mathscr{G}$ in a Cartan geometry $(\mathscr{G},\omega)$ of type $(G,H)$, the \emph{development} $\gamma_G$ of $\gamma$ is the unique (piecewise smooth) path $\gamma_G\colon [0,1]\to G$ such that $\gamma_G(0)=e$ and $\gamma^*\omega=\gamma_G^*\omega_G$.
\end{Definition}

Development allows us to identify paths in a Cartan geometry with paths in the model, and if we fix a pretend ``identity element'' $\mathscr{e}\in\mathscr{G}$, we can even give a kind of correspondence between elements of $\mathscr{G}$ and elements of $G$.

\begin{Definition}For a Cartan geometry $(\mathscr{G},\omega)$ of type $(G,H)$ and points $\mathscr{e},\mathscr{g}\in\mathscr{G}$, we say that $g\in G$ is a \emph{development of $\mathscr{g}$ from $\mathscr{e}$} if and only if there exists a path $\gamma\colon [0,1]\to\mathscr{G}$ and $h\in H$ such that $\gamma(0)=\mathscr{e}$, $\gamma(1)h=\mathscr{g}$, and $\gamma_G(1)h=g$.\end{Definition}

Developments of elements in $\mathscr{G}$ are usually not unique. Thankfully, the \emph{holonomy group} $\Hol_\mathscr{e}(\mathscr{G},\omega)$ tells us precisely how this happens: if $g$ is a development of $\mathscr{g}$ from $\mathscr{e}$, then the set of all possible developments from $\mathscr{e}$ to $\mathscr{g}$ is precisely $\Hol_\mathscr{e}(\mathscr{G},\omega)g$.

\begin{Definition}For a Cartan geometry $(\mathscr{G},\omega)$ of type $(G,H)$, the \emph{holonomy group} of $(\mathscr{G},\omega)$ at $\mathscr{e}\in\mathscr{G}$ is the subgroup of $G$ given by \[\Hol_\mathscr{e}(\mathscr{G},\omega):=\left\{\gamma_G(1)h_\gamma^{-1}\in G \,\Big|\, \begin{array}{@{}l@{}}\gamma\colon [0,1]\to\mathscr{G}\text{ is a path such that,} \\ \text{for }h_\gamma\in H,\,\gamma(0)=\mathscr{e}=\gamma(1)h_\gamma^{-1}\end{array}\right\}.\]
\end{Definition}

Because $\Hol_\mathscr{e}(\mathscr{G},\omega)$ completely encapsulates the ambiguity in taking developments from $\mathscr{e}$, we can assign unique developments to elements of $\mathscr{G}$ when $\Hol_\mathscr{e}(\mathscr{G},\omega)=\{e\}$. In other words, whenever $\Hol_\mathscr{e}(\mathscr{G},\omega)=\{e\}$, we get a geometric map $\mathrm{dev}_\mathscr{e}\colon (\mathscr{G},\omega)\to (G,\MC{G})$, which we could reasonably call the \emph{developing map based at $\mathscr{e}$}, given by $\gamma(1)h\mapsto\gamma_G(1)h$ for (piecewise smooth) paths $\gamma \colon [0,1]\to\mathscr{G}$ starting at $\gamma(0)=\mathscr{e}$ and $h\in H$. When the geometry is flat and the base manifold is simply connected, the holonomy group is always trivial, and the induced map of $\mathrm{dev}_\mathscr{e}$ on the base manifold is precisely the developing map in the usual sense for locally homogeneous geometric structures.

Crucially, it follows that the curvature vanishes over an open set if and only if loops that are null-homotopic in that open set have trivial holonomy.

\section{Some useful tools for working with holonomy}\label{devtools}
To make appropriate use of holonomy and development for our present purposes, we will need a way to tell when a sequence of holonomies converges, as well as a way to reverse the process of development for sufficiently small paths. In this section, we will provide a couple of results that will allow us to do this, which are useful in their own right and will almost certainly see significant use in subsequent works.

\subsection{Convergence of developments}
For a Cartan geometry $(\mathscr{G},\omega)$ of type $(G,H)$, an inner product $\mathrm{g}$ on the Lie algebra $\mathfrak{g}$ determines a Riemannian metric $\mathrm{g}_\omega$ on $\mathscr{G}$ given by $\mathrm{g}_\omega(\xi_1,\xi_2):=\mathrm{g}(\omega(\xi_1),\omega(\xi_2))$; we call such a Riemannian metric on $\mathscr{G}$ an \emph{$\omega$-constant Riemannian metric}. Given an element $\mathscr{e}\in\mathscr{G}$, we can use an arbitrary $\omega$-constant Riemannian metric $\mathrm{g}_\omega$ to imbue the set $\Path_\mathscr{e}(\mathscr{G})$ of (piecewise smooth) paths $\gamma \colon [0,1]\to\mathscr{G}$ in $\mathscr{G}$ starting at $\gamma(0)=\mathscr{e}$ with a useful topology.

\begin{Definition}The \emph{first-order compact-open topology} on \[\Path_\mathscr{e}(\mathscr{G}):=\bigl\{\gamma\in C^0([0,1];\mathscr{G})\mid  \gamma\text{ is a piecewise smooth path and }\gamma(0)=\mathscr{e} \bigr\},\] the space of piecewise smooth paths $\gamma \colon [0,1]\to\mathscr{G}$ starting at $\gamma(0)=\mathscr{e}$, is the topology induced by the metric \[d_{\mathrm{g}_\omega}(\gamma_1,\gamma_2):=\sup_{t\in[0,1]} (\mathrm{dist}_{\mathrm{g}_\omega}(\gamma_1(t), \gamma_2(t)) )+\sup_{t\in[0,1]}\|\omega(\dot{\gamma}_1(t))-\omega(\dot{\gamma}_2(t))\|_\mathrm{g},\] where the second supremum is taken over all $t\in[0,1]$ such that both $\dot{\gamma}_1(t)$ and $\dot{\gamma}_2(t)$ are well defined.
\end{Definition}

This is essentially just an extension of the usual $C^1$ compact-open topology on smooth paths so that it is defined on piecewise smooth paths instead.

When we say that a sequence of paths converges in $\Path_\mathscr{e}(\mathscr{G})$, we will mean that the sequence converges with respect to the first-order compact-open topology.

\begin{Lemma}[continuous dependence on parameters]\label{cdop} Suppose $(f_\lambda)$ is a net of functions from $[0,1]$ to $\mathfrak{g}$ that converges uniformly to some function $f_\infty\colon [0,1]\to\mathfrak{g}$. If there exists a solution $\gamma_\infty\in\Path_\mathscr{e}(\mathscr{G})$ to the ODE given by $\omega(\dot{\gamma}_\infty)=f_\infty$, and solutions $\gamma_\lambda\in\Path_\mathscr{e}(\mathscr{G})$ to $\omega(\dot{\gamma}_\lambda)=f_\lambda$ for each $\lambda$, then $(\gamma_\lambda)$ converges to $\gamma_\infty$ in $\Path_\mathscr{e}(\mathscr{G})$.\end{Lemma}
\begin{proof}As the labeling implies, this result is effectively just a version of the usual ``continuous dependence on parameters'' lemma that one typically encounters for ODEs, except that it is adapted to our current setting. The main tool we will need is the ``fundamental lemma''~1.5.1 from \cite[Chapter~2]{HCartan1971}, which tells us that, with the norm $\|\cdot\|$ on a Banach space $E$, a subset $U\subseteq E$, and an interval $I$, if a function $f\colon I\times U\to E$ is $k$-Lipschitz, $\gamma \colon I\to U$ solves the ODE given by $\dot{\gamma}(t)=f(t,\gamma(t))$, and $\bar{\gamma}\colon I\to U$ satisfies $\|\dot{\bar{\gamma}}(t)-f(t,\bar{\gamma}(t))\|\leq\varepsilon$ for all $t\in I$, then for $\tau\in I$, the distance $\|\gamma(t)-\bar{\gamma}(t)\|$ satisfies the bound \[\|\gamma(t)-\bar{\gamma}(t)\|\leq\|\gamma(\tau)-\bar{\gamma}(\tau)\|{\rm e}^{k|t-\tau|}+\varepsilon\frac{{\rm e}^{k|t-\tau|}-1}{k}.\]

The strategy, in essence, is just to apply this ``fundamental lemma'' in finitely many coordinate balls inside charts along $\gamma_\infty$. This will show that $\gamma_\lambda$ converges uniformly to $\gamma_\infty$ in each of the finitely many balls, hence it converges to $\gamma_\infty$ uniformly, and since the derivatives converge uniformly by assumption, it follows that $\gamma_\lambda$ converges to $\gamma_\infty$ in the first-order compact-open topology.

To be more explicit, let us cover $\gamma_\infty([0,1])$ in finitely many coordinate balls $B_0,\dots,B_n$, and cover $[0,1]$ by closed subintervals $I_0,\dots,I_n$ with $\max(I_i)\in I_{i+1}$ and $\gamma_\infty(I_i)\subseteq B_i$ for each $i$. Since $f_\infty(t)$ is constant in $\mathscr{g}\in\mathscr{G}$, we can find $k_i>0$ for which the coordinate expressions for the smooth vector fields $\omega^{-1}(f_\infty(t))$ are $k_i$-Lipschitz on $\bar{B}_i$ with respect to the Euclidean norm $\|\cdot\|_\text{coord}$ induced by the coordinates; for clarity, note that this norm only gives a metric on the chart, not all of $\mathscr{G}$, but since $\bar{B}_i$ is compact, the restriction of the coordinate metric to $\bar{B}_i$ is bi-Lipschitz equivalent to the restriction of an $\omega$-constant metric, % since the coordinate norms of the omega-constant unit vectors are bounded above and below on the balls,
so uniform convergence on $\bar{B}_i$ with respect to this metric implies uniform convergence with respect to an $\omega$-constant one and vice versa. By compactness of $[0,1]$ and $\bar{B}_i$,
\[\varepsilon_i(\lambda):=\max_{\substack{t\in[0,1] \\ \mathscr{g}\in\bar{B}_i}}\bigl\|\omega_\mathscr{g}^{-1}(f_\lambda(t))-\omega_\mathscr{g}^{-1}(f_\infty(t))\bigr\|_\text{coord},\]
exists for each $i$ and $\lambda$, and it converges to $0$ as $\lambda\rightarrow\infty$ because $f_\lambda$ converges uniformly to $f_\infty$ by assumption. Starting with $t_0:=0$, so that $\|\gamma_\lambda(t_0)-\gamma_\infty(t_0)\|_\text{coord}=0$, the ``fundamental lemma'' of \cite{HCartan1971} tells us that, for $\lambda$ sufficiently far along toward $\infty$,
\[\|\gamma_\lambda(t)-\gamma_\infty(t)\|_\text{coord}\leq\varepsilon_0(\lambda)\frac{{\rm e}^{k_0t}-1}{k_0}\leq\varepsilon_0(\lambda)\frac{{\rm e}^{k_0}-1}{k_0}\]
for all $t\in I_0$, so as $\lambda\rightarrow\infty$, $\gamma_\lambda|_{I_0}$ converges uniformly to $\gamma_\infty|_{I_0}$. Next, choose $t_1\in I_0\cap I_1$; by taking $\lambda$ to be sufficiently far along toward $\infty$, we have just seen (now using the coordinate metric in $B_1$ rather than $B_0$) that we can get $\|\gamma_\lambda(t_1)-\gamma_\infty(t_1)\|_\text{coord}$ as close to $0$ as we would like, so using the ``fundamental lemma'' of \cite{HCartan1971}, we get that
\[\|\gamma_\lambda(t)-\gamma_\infty(t)\|_\text{coord}\leq\|\gamma_\lambda(t_1)-\gamma_\infty(t_1)\|_\text{coord}{\rm e}^{k_1}+\varepsilon_1(\lambda)\frac{{\rm e}^{k_1}-1}{k_1}\]
for all $t\in I_1$, hence $\gamma_\lambda|_{I_1}$ converges uniformly to $\gamma_\infty|_{I_1}$ as $\lambda\rightarrow\infty$. Iteratively, we repeat this process, picking $t_i\in I_{i-1}\cap I_i$ and then using the uniform convergence over $I_{i-1}$ and the ``fundamental lemma'' from \cite{HCartan1971} to show that, for $\lambda$ sufficiently far along toward $\infty$, we get
\[\|\gamma_\lambda(t)-\gamma_\infty(t)\|_\text{coord}\leq\|\gamma_\lambda(t_i)-\gamma_\infty(t_i)\|_\text{coord}\,{\rm e}^{k_i}+\varepsilon_i(\lambda)\frac{{\rm e}^{k_i}-1}{k_i}\]
for all $t\in I_i$, hence $\gamma_\lambda|_{I_i}$ converges uniformly to $\gamma_\infty|_{I_i}$ for each of the finitely many $i$. As noted above, this implies that $\gamma_\lambda$ converges uniformly to $\gamma_\infty$ on all of $[0,1]$, and since, by assumption, the derivatives converge uniformly as well, $\gamma_\lambda$ converges to $\gamma_\infty$ in $\Path_\mathscr{e}(\mathscr{G})$.\end{proof}

This naturally leads to the following corollary, which seems implicit in much of the existing literature.

\begin{Proposition}\label{devcont} The development map $(\cdot)_G\colon \Path_\mathscr{e}(\mathscr{G})\to\Path_e(G)$ given by $\gamma\mapsto\gamma_G$ is continuous.\end{Proposition}
\begin{proof}Since $\Path_\mathscr{e}(\mathscr{G})$ is a metric space, it suffices to prove that if $\gamma_k\rightarrow\gamma_\infty$ in $\Path_\mathscr{e}(\mathscr{G})$, then $(\gamma_k)_G\rightarrow(\gamma_\infty)_G$ in $\Path_e(G)$. But $\gamma_k\rightarrow\gamma_\infty$ in $\Path_\mathscr{e}(\mathscr{G})$ implies that $\omega(\dot{\gamma}_k)=\MC{G}((\dot{\gamma}_k)_G)$ converges uniformly to $\omega(\dot{\gamma}_\infty)=\MC{G}((\dot{\gamma}_\infty)_G)$, so by Lemma~\ref{cdop}, $(\gamma_k)_G\rightarrow(\gamma_\infty)_G$ in $\Path_e(G)$. \end{proof}

In particular, it follows from Proposition~\ref{devcont} that a sequence of paths converging to $\mathscr{e}$ in $\Path_\mathscr{e}(\mathscr{G})$ must have developments converging to the identity in $G$.

\begin{Lemma}\label{shrinkingholonomy}Suppose $(\gamma_k)$ is a sequence in $\Path_\mathscr{e}(\mathscr{G})$. If $(\gamma_k)$ converges in $\Path_\mathscr{e}(\mathscr{G})$ to the constant path at $\mathscr{e}$, then the developments $(\gamma_k)_G(1)$ converge to the identity element $e\in G$.\end{Lemma}

For the purposes of using this lemma, it will be helpful to keep two particular facts in mind.

First, the endpoint of the development of a path does not depend on parametrization: if $\tau$ is a (piecewise smooth) reparametrization of $[0,1]$ with $\tau(0)=0$ and $\tau(1)=1$, then $\gamma^*\omega=\gamma_G^*\omega_G$ implies \[(\gamma\circ\tau)^*\omega=\tau^*\gamma^*\omega=\tau^*\gamma_G^*\omega_G=(\gamma_G\circ\tau)^*\omega_G,\] so $(\gamma\circ\tau)_G=\gamma_G\circ\tau$ and $(\gamma\circ\tau)_G(1)=\gamma_G(\tau(1))=\gamma_G(1)$. Therefore, even if a sequence of paths in $\Path_\mathscr{e}(\mathscr{G})$ converges to the constant path \textit{modulo parametrization}, the endpoints of their developments will still converge to the identity in $G$.

Second, the lemma is applicable even when the converging paths we are considering are on the base manifold: if $(\gamma_k)$ is a sequence of paths in $\mathscr{G}$ such that $\gamma_k(0)\in\mathscr{e}H$ and the paths~$q_{{}_H}(\gamma_k)$ in the base manifold converge to the constant path at $q_{{}_H}(\mathscr{e})$ to first order, then we can find paths $\beta_k\colon [0,1]\to H$ such that $\gamma_k(0)\beta_k(0)^{-1}=\mathscr{e}$ and $\gamma_k\beta_k^{-1}\rightarrow\mathscr{e}$ in $\Path_\mathscr{e}(\mathscr{G})$. In this situation, we should take particular note of how the sequence $(\beta_k)$ of paths in $H$ plays a similar role to what we call \textit{ballast sequences}---the original term ``holonomy sequence'', as used in \cite{CapMelnick2013,Frances2007,MelnickNeusser2016}, among several other places, is clearly inappropriate given the present context, since they are not sequences of holonomies for the geometry---in Section~\ref{snares} below: while the paths $\gamma_k$ might not converge to anything, we can use the paths $\beta_k$ to counterbalance them so that the end result~$\gamma_k\beta_k^{-1}$ still converges.

\subsection{Existence of antidevelopments}
It would be convenient to be able to reverse the process of development: given a piecewise smooth path $\gamma \colon [0,1]\to G$ with $\gamma(0)=e$, we would like to be able to construct a corresponding path $\gamma_\mathscr{G}$ in $\mathscr{G}$ such that $\gamma$ is the development of $\gamma_\mathscr{G}$. We can reasonably call such a path---if it exists---an \emph{antidevelopment} of $\gamma$.

\begin{Definition}For a Cartan geometry $(\mathscr{G},\omega)$ of type $(G,H)$ and a path $\gamma\in\Path_e(G)$, an \emph{antidevelopment} for $\gamma$ starting at $\mathscr{e}\in\mathscr{G}$ is a path $\gamma_\mathscr{G}\in\Path_\mathscr{e}(\mathscr{G})$ such that $(\gamma_\mathscr{G})_G=\gamma$.\end{Definition}

Unlike with developments, antidevelopments do not always exist. Conceptually, we can regard this asymmetry as a consequence of the fact that $\omega$-constant metrics $\mathrm{g}_\omega$ for a Cartan geometry $(\mathscr{G},\omega)$ are not necessarily complete, whereas $\omega_G$-constant metrics on $G$ are always complete. In other words, there will often be paths that ``escape'' $\mathscr{G}$ despite having finite length with respect to some $\mathrm{g}_\omega$, but this can never happen for finite-length paths in $G$ because $\omega_G$-constant metrics on $G$ must be complete.

However, if we fix a particular element $\mathscr{e}\in\mathscr{G}$ and an $\omega$-constant metric $\mathrm{g}_\omega$, then we can always find a sufficiently small radius $r_\mathscr{e}>0$ for which the metric ball $B_{r_\mathscr{e}}(\mathscr{e})$ has compact closure. This radius will usually depend on the choice of element $\mathscr{e}$ unless $\mathrm{g}_\omega$ is complete, but it will always be greater than $0$. Paths starting at $\mathscr{e}$ of length less than $r_\mathscr{e}$ will be trapped inside the compact subspace $\overline{B_{r_\mathscr{e}}(\mathscr{e})}$, so if $\gamma \colon [0,1]\to G$ is a path such that $\gamma(0)=e$ whose length with respect to $\mathrm{g}_{\omega_G}$ is less than $r_\mathscr{e}$, then its antidevelopment $\gamma_\mathscr{G}$ will not be able to escape $\overline{B_{r_\mathscr{e}}(\mathscr{e})}$. In particular, we have the following.

\begin{Lemma}\label{shortantidevs} Suppose $(\mathscr{G},\omega)$ is a Cartan geometry of type $(G,H)$ and $\mathrm{g}$ is an inner product on~$\mathfrak{g}$. For each $\mathscr{e}\in\mathscr{G}$, there exists an $r_\mathscr{e}>0$ such that every path $\gamma\in\Path_e(G)$ whose length with respect to $\mathrm{g}_{\omega_G}$ is less than $r_\mathscr{e}$ has a well-defined antidevelopment $\gamma_\mathscr{G}$ starting at~$\mathscr{e}$.\end{Lemma}

In summary, sufficiently short paths $\gamma\in\Path_e(G)$ will always have antidevelopments at a~given point $\mathscr{e}\in\mathscr{G}$.

\section{Isotropies that ``ensnare'' a region of the model}\label{snares}
In \cite{CapMelnick2013,MelnickNeusser2016}, the overall strategy for demonstrating local flatness near a given higher-order fixed point was to find an open region of the geometry whose points were pulled toward that higher-order fixed point in a particular way. Here, we will establish two criteria to capture the behavior underlying these results when the higher-order fixed point is \textit{isolated}, meaning that some neighborhood of it contains no other higher-order fixed points with the same isotropy up to conjugation. These criteria---codified in Definition \ref{ensnaredef}---are entirely determined by the behavior of the isotropy in the relevant model geometry, which we will not require to be parabolic in general. Because this behavior only depends on the isotropy, this allows us to transfer this desired behavior to every Cartan geometry admitting an automorphism with such isotropy.

\subsection{Criteria for the isotropy}
\begin{figure}[!ht]
\centering
\includegraphics[width=0.6\textwidth]{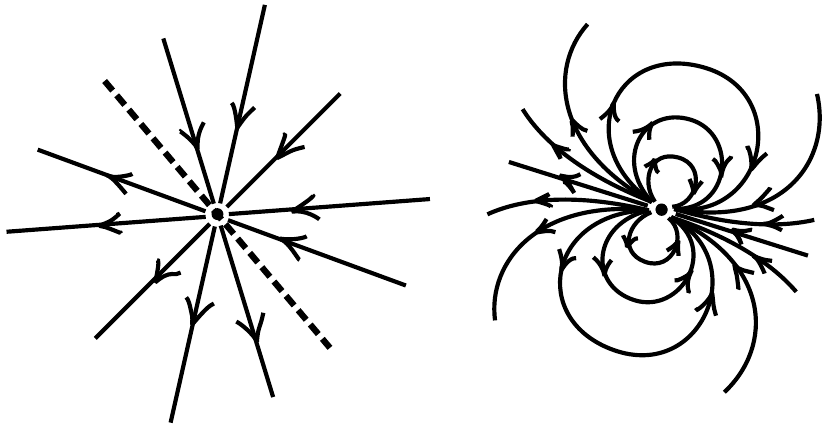}
\caption{Sketches of automorphisms with an isolated higher-order fixed point for 2-dimensional real projective geometry $(\PGL_3\mathbb{R},P)$ (left) and 1-dimensional complex projective geometry $(\PGL_2\mathbb{C},P)$ (right).}
\label{projsnareimg}
\end{figure}

To understand the expected behavior of automorphisms near an isolated higher-order fixed point, which we will take as inspiration for the ideas below, it is useful to look at some examples on low-dimensional model geometries. We have sketched such examples for real and complex projective geometry in Figure~\ref{projsnareimg}.

Given an automorphism $\alpha$ with a higher-order fixed point, the results in \cite{CapMelnick2013,MelnickNeusser2016} suggest we should look for particular types of \emph{ballast sequences} to find an open set on which the curvature vanishes.

\begin{Definition}Let $\alpha\in\Aut(\mathscr{G},\omega)$ be an automorphism for a Cartan geometry $(\mathscr{G},\omega)$ of type $(G,H)$. A sequence $(b_k)$ in $H$ is called a \emph{ballast sequence for $\alpha$ at $\mathscr{g}\in\mathscr{G}$ with attractor $\mathscr{e}\in\mathscr{G}$} if and only if there exists a sequence $(\mathscr{g}_k)$ in $\mathscr{G}$ such that $\mathscr{g}_k\rightarrow\mathscr{g}$ and $\alpha^k(\mathscr{g}_k)b_k^{-1}\rightarrow\mathscr{e}$.\end{Definition}

In particular, for isolated higher-order fixed points, the techniques in \cite{CapMelnick2013,MelnickNeusser2016} suggest looking for ballast sequences $(b_k)$ contained in the Levi subgroup $G_0<P$ of our parabolic isotropy~$P$ such that the restriction $\Ad_{b_k}|_{\mathfrak{g}_-}$ of the adjoint action to $\mathfrak{g}_-$ converges to~$0$. We should think of this as an assumption on the behavior of the derivatives of the automorphism on the tangent bundle over the base manifold $M$. To understand why, suppose that $(b_k)$ is such a ballast sequence at~$\mathscr{g}$, so that $\alpha^k(\mathscr{g}_k)b_k^{-1}\rightarrow\mathscr{e}$ for some sequence $(\mathscr{g}_k)$ in $\mathscr{G}$ converging to $\mathscr{g}$ and $\Ad_{b_k}(v)\rightarrow 0$ for all $v\in\mathfrak{g}_-$. For sufficiently small $v\in\mathfrak{g}_-$, we therefore get a well-defined sequence of elements $\exp\bigl(\omega^{-1}(v)\bigr)\mathscr{g}_k$ converging to $\exp\bigl(\omega^{-1}(v)\bigr)\mathscr{g}$, where $\exp\bigl(\omega^{-1}(v)\bigr)$ denotes the time-1 flow along the $\omega$-constant vector field $\omega^{-1}(v)$, so that
\[\alpha^k\bigl(\exp\bigl(\omega^{-1}(v)\bigr)\mathscr{g}_k\bigr)b_k^{-1} =\exp\bigl(\omega^{-1}(\Ad_{b_k}v)\bigr)\bigl(\alpha^k(\mathscr{g}_k)b_k^{-1}\bigr)\rightarrow\mathscr{e}.\]
Taking such $v$ to be in a bounded neighborhood $V$ of $0\in\mathfrak{g}_-$, we can make $\Ad_{b_k}(V)$ as close to $0$ as we would like for $k$ sufficiently large, so for every neighborhood of $q_{{}_P}(\mathscr{e})$, we may find a large enough $k$ that $\alpha^k q_{{}_P}\bigl(\exp\bigl(\omega^{-1}(V)\bigr)\mathscr{g}_k\bigr)$ is contained in that neighborhood. In other words, the existence of such a ballast sequence for each $\mathscr{g}\in q_{{}_P}^{-1}(U)$ over some open set $U\subseteq M$ guarantees that the iterates $\alpha^k$, viewed as maps on the base manifold, converge \textit{locally uniformly} to the constant map at $q_{{}_P}(\mathscr{e})$ on $U$. Moreover, for whatever $w\in\mathfrak{g}_-$ we choose, we can define $\exp\bigl(t\omega^{-1}(w)\bigr)\exp\bigl(\omega^{-1}(v)\bigr)\mathscr{g}_k$ for sufficiently small $t\in\mathbb{R}$, and
\begin{gather*}
\alpha^k\bigl(\exp\bigl(t\omega^{-1}(w)\bigr)\exp\bigl(\omega^{-1}(v)\bigr)\mathscr{g}_k\bigr)b_k^{-1} \\
\qquad =\exp\bigl(t\omega^{-1}(\Ad_{b_k}w)\bigr)\exp\bigl(\omega^{-1}(\Ad_{b_k}v)\bigr) \bigl(\alpha^k(\mathscr{g}_k)b_k^{-1}\bigr)\rightarrow\mathscr{e};
\end{gather*}
restricting our possible choices of $w$ to a bounded set in $\mathfrak{g}_-$, we can again get all $\Ad_{b_k}(w)$ as close to $0$ as we would like, so the corresponding curves $t\mapsto \alpha^k\bigl(q_{{}_P}\bigl(\exp\bigl(t\omega^{-1}(w)\bigr)\exp\bigl(\omega^{-1}(v)\bigr)\mathscr{g}_k\bigr)\bigr)$ converge to first order to the constant curve at $q_{{}_P}(\mathscr{e})$. Thus, we are looking for automorphisms~$\alpha$ with a~region~$U$ over which iterates of the derivatives of~$\alpha$ converge locally uniformly to the zero vector over some point of the base manifold.

Another necessary component in the techniques of~\cite{CapMelnick2013,MelnickNeusser2016} is the presence of certain paths that get shrunk to a point of the base manifold by the automorphism. The flat regions constructed in those papers are covered by such shrinking paths, so in particular, there should be at least one of them. We will only need to specify one such path, which we will call~$\zeta_U$, in the base manifold of the model geometry, and moreover, we will not need it to be the quotient image of an exponential path, though (for convenience) we will choose for it to be so in all of the applications in Section~\ref{applications}.

Combining these two criteria leads us to the following definition.

\begin{Definition}\label{ensnaredef} We say that $a\in H$ \emph{ensnares} an $a$-invariant connected open subset $U\subseteq G/H$ if and only if the following two criteria hold:
\begin{itemize}\itemsep=0pt
\item Over $U$, the derivative maps $a^k_*\colon T(G/H)\to T(G/H)$ of the iterates $a^k$ converge locally uniformly to the constant map at the zero vector $0_{q_{{}_H}(e)}$ at $q_{{}_H}(e)$ as $k\rightarrow+\infty$.
\item There exists a path $\zeta_U\colon [0,1]\to G/H$ such that $\zeta_U(0)=q_{{}_H}(e)$, $\zeta_U(1)\in U$, and, for some sequence $(\tau_k)$ of orientation-preserving reparametrizations of $[0,1]$, $a^k(\zeta_U\circ\tau_k)$ converges to the constant path at $q_{{}_H}(e)$ to first order as $k\rightarrow+\infty$.
\end{itemize}\end{Definition}

Assuming local uniform convergence over $U$ of the derivatives of the iterates of $a$ on $G/H$ provides a convenient way to guarantee that paths inside $U$ are pulled to $q_{{}_H}(e)$. More precisely, it tells us that whenever we have a path in $U$, its iterates under $a$ must converge to the constant path at $q_{{}_H}(e)$ to first order.

\begin{figure}[!ht]
\centering
\includegraphics[width=0.7\textwidth]{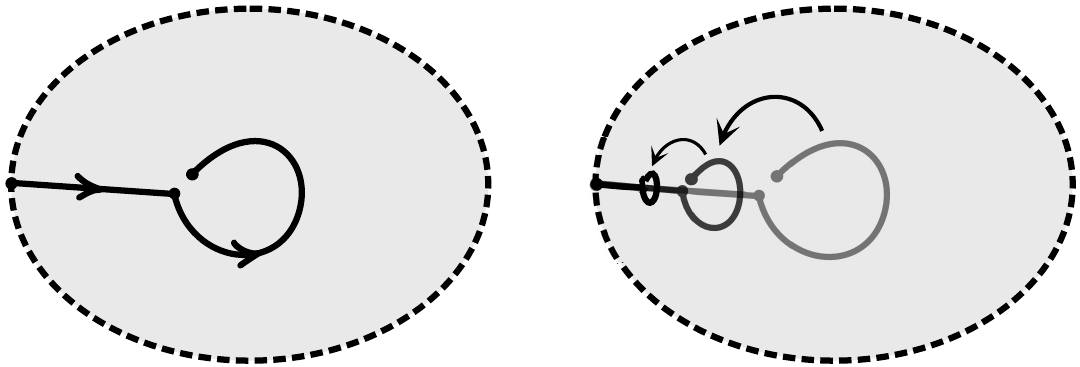}
\put(-92,82){\makebox(0,0)[lb]{\small $a$}}
\put(-120,69){\makebox(0,0)[lb]{\small $a$}}
\put(-228,37){\makebox(0,0)[lb]{\small $q_{{}_H}(\delta)$}}
\put(-298,39){\makebox(0,0)[lb]{\small $\zeta_U$}}
\put(-228,90){\makebox(0,0)[lb]{\small $U$}}

\caption{A concatenation of the form $\zeta_U\star q_{{}_H}(\delta)$ (left) and some of its iterates under $a$ (right).}
\label{snareimg}
\end{figure}
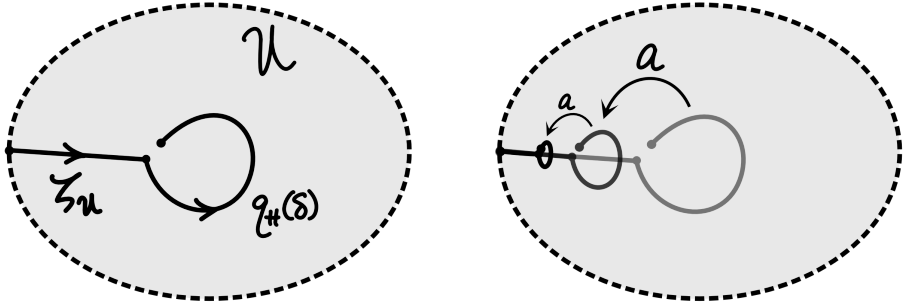

For paths $\gamma_1$ and $\gamma_2$ such that $\gamma_1(1)=\gamma_2(0)$, so that the endpoint of the first is the starting point of the second, we denote by $\gamma_1\star\gamma_2$ their \emph{concatenation}
\[\gamma_1\star\gamma_2(t):= \begin{cases}\gamma_1(2t) & \text{if }t\in\bigl[0,\tfrac{1}{2}\bigr], \\
\gamma_2(2t-1) & \text{if }t\in\bigl[\tfrac{1}{2},1\bigr].\end{cases}
\]
From what Definition \ref{ensnaredef} tells us about $\zeta_U$, together with the aforementioned convergence of paths in $U$ under forward iterates of $a$, it follows that concatenations of the form $\zeta_U\star q_{{}_H}(\delta)$, where $\delta\colon [0,1]\to q_{{}_H}^{-1}(U)$ is a path such that $q_{{}_H}(\delta(0))=\zeta_U(1)$, get pulled into $q_{{}_H}(e)$, so that $a^k((\zeta_U\circ\tau_k)\star q_{{}_H}(\delta))$ also converges to $q_{{}_H}(e)$ to first order. As illustrated in Figure~\ref{snareimg}, the resulting picture can end up looking a bit like a snare being pulled closed, hence the name.

Before we proceed, it is worth noting again that, while the two criteria of Definition~\ref{ensnaredef} were designed to generalize the behavior used by \cite{CapMelnick2013,MelnickNeusser2016} for their results on isolated higher-order fixed points, this parabolic context is not actually required in the definition. In particular, even though our main applications in this paper happen to involve automorphisms of parabolic geometries with isolated higher-order fixed points, we will not need to use this as an explicit assumption for the general results below. For example, in affine geometry $(\mathrm{Aff}(m),\GLin_m\mathbb{R})$, the contracting homothety $a=\lambda\mathds{1}\in\GLin_m\mathbb{R}$ (where $|\lambda|<1$) ensnares the whole model space $U=\mathrm{Aff}(m)/\GLin_m\mathbb{R}\cong\mathbb{R}^m$, with $\zeta_U$ the constant path at $q_{{}_{\GLin_m\mathbb{R}}}(e)$.

\subsection{Existence of ensnared region with trivial holonomy} Since we can pull in concatenations of the form $\zeta_U\star q_{{}_H}(\delta)$ to be as short as we want, we can show that they always have well-defined antidevelopments when there exists an automorphism with ensnaring isotropy.

\begin{Lemma}\label{snarelemma} If $(\mathscr{G},\omega)$ is a Cartan geometry of type $(G,H)$, $\mathscr{e}\in\mathscr{G}$, and $\alpha\in\Aut(\mathscr{G},\omega)$ such that $\alpha(\mathscr{e})=\mathscr{e}a$ for some $a\in H$ that ensnares $U\subseteq G/H$, then every path of the form $\hat{\zeta}_U\star\delta$, where $\hat{\zeta}_U\in\Path_e(G)$ is a lift of $\zeta_U$ and $\delta\colon [0,1]\to q_{{}_H}^{-1}(U)$ such that $\delta(0)=\hat{\zeta}_U(1)$, has a~well-defined antidevelopment \smash{$\bigl(\hat{\zeta}_U\star\delta\bigr)_\mathscr{G}$} starting at $\mathscr{e}$.\end{Lemma}
\begin{proof}Since \smash{$a^k((\zeta_U\circ\tau_k)\star q_{{}_H}(\delta))=a^k q_{{}_H}\bigl(\bigl(\hat{\zeta}_U\circ\tau_k\bigr)\!\star\delta\bigr)$} converges to the constant path at $q_{{}_H}(e)$ to first order, there exist $\beta_k\colon [0,1]\to H$ such that $\beta_k(0)=a^k$ and \smash{$a^k\bigl(\hat{\zeta}_U\star\delta\bigr)\beta_k^{-1}$} converges, modulo parametrization, to the constant path at $e$ in $\Path_e(G)$. As a result, for each $r>0$, there is some $k\geq 0$ for which $a^k\bigl(\hat{\zeta}_U\star\delta\bigr)\beta_k^{-1}$ has length less than~$r$.

By Lemma~\ref{shortantidevs}, there exists an $r_\mathscr{e}>0$ for which all paths in $\Path_e(G)$ of length less than $r_\mathscr{e}$ have well-defined antidevelopments starting at $\mathscr{e}\in\mathscr{G}$. Thus, $a^k\bigl(\hat{\zeta}_U\star\delta\bigr)\beta_k^{-1}$ has an antidevelopment $\bigl(a^k\bigl(\hat{\zeta}_U\star\delta\bigr)\beta_k^{-1}\bigr)_\mathscr{G}$ starting at $\mathscr{e}$ for all sufficiently large $k$, hence \[\bigl(\hat{\zeta}_U\star\delta\bigr)_\mathscr{G}:= \alpha^{-k}\bigl(a^k\bigl(\hat{\zeta}_U\star\delta\bigr)\beta_k^{-1}\bigr)_\mathscr{G}\beta_k\] is a~well-defined path starting at $\alpha^{-k}(\mathscr{e})\beta_k(0)=\mathscr{e}a^{-k}a^k=\mathscr{e}$ with development \begin{align*}\bigl(\bigl(\hat{\zeta}_U\star\delta\bigr)_\mathscr{G}\bigr)_G & =\bigl(\alpha^{-k}\bigl(a^k\bigl(\hat{\zeta}_U\star\delta\bigr)\beta_k^{-1}\bigr)_\mathscr{G}\beta_k\bigr)_G
   =\beta_k(0)^{-1}\bigl(\bigl(a^k\bigl(\hat{\zeta}_U\star\delta\bigr)\beta_k^{-1}\bigr)_\mathscr{G}\bigr)_G\beta_k \\
 & =a^{-k}\bigl(a^k\bigl(\hat{\zeta}_U\star\delta\bigr)\beta_k^{-1}\bigr)\beta_k   =\hat{\zeta}_U\star\delta.\tag*{\qed}
 \end{align*}\renewcommand{\qed}{}
 \end{proof}

Whenever $(\mathscr{G},\omega)$ has an automorphism $\alpha\in\Aut(\mathscr{G},\omega)$ with isotropy $a\in H$ that ensnares $U\subseteq G/H$, we can use Lemma~\ref{snarelemma} to construct an $H$-invariant subspace \[q_{{}_H}^{-1}(U)_\mathscr{G}:=\bigl\{\bigl(\hat{\zeta}_U\star\delta\bigr)_\mathscr{G}(1)h\mid \delta\in\Path_{\hat{\zeta}_U(1)}\bigl(q_{{}_H}^{-1}(U)\bigr),\,h\in H\bigl\}\subseteq\mathscr{G}\] given by the endpoints of antidevelopments for concatenations of the form $\hat{\zeta}_U\star\delta$ together with their right-translations by elements of $H$. Note that $q_{{}_H}^{-1}(U)_\mathscr{G}$ is necessarily an open subset of~$\mathscr{G}$: for each $\delta(1)\in q_{{}_H}^{-1}(U)$, sufficiently short paths starting at $\delta(1)$ stay in the open set $q_{{}_H}^{-1}(U)$, so each \smash{$\bigl(\hat{\zeta}_U\star\delta\bigr)_\mathscr{G}(1)\in q_{{}_H}^{-1}(U)_\mathscr{G}$} is the center of a small open ball contained in~$q_{{}_H}^{-1}(U)_\mathscr{G}$.

We should think of $q_{{}_H}^{-1}(U)_\mathscr{G}$ as the analogue in $\mathscr{G}$ for the open subset $q_{{}_H}^{-1}(U)$ in~$G$, and we~should hope for it to exhibit similar behavior with respect to the automorphism $\alpha\in\Aut(\mathscr{G},\omega)$ to what $q_{{}_H}^{-1}(U)$ does for the ensnaring isotropy $a\in H$. Indeed, as we mentioned in the proof of the lemma, each path of the form \smash{$\hat{\zeta}_U\star\delta$} admits a sequence of paths $\beta_k\colon [0,1]\to H$ with $\beta_k(0)=a^k$ such that \smash{$a^k\bigl(\hat{\zeta}_U\star\delta\bigr)\beta_k^{-1}$} converges to the constant path $e$ in $\Path_e(G)$ modulo parametrization, and \[\bigl(\alpha^k\bigl(\hat{\zeta}_U\star\delta\bigr)_\mathscr{G}\beta_k^{-1}\bigr)_G =\beta_k(0)\bigl(\hat{\zeta}_U\star\delta\bigr)\beta_k^{-1} =a^k\bigl(\hat{\zeta}_U\star\delta\bigr)\beta_k^{-1},
\] so we must also have that \smash{$\alpha^k\bigl(\hat{\zeta}_U\star\delta\bigr)_\mathscr{G}\beta_k^{-1}$} converges to $\mathscr{e}$ in $\Path_\mathscr{e}(\mathscr{G})$ modulo parametrization. Consequently, the derivatives $\alpha^k_*$ on the base manifold must converge locally uniformly to $0_{q_{{}_H}(\mathscr{e})}$ over the connected open set \smash{$U_\mathscr{G}:=q_{{}_H}\bigl(q_{{}_H}^{-1}(U)_\mathscr{G}\bigr)$}, just like in the model case. Using this, we can prove the key result for this section: when $a^kga^{-k}\rightarrow e$ implies $g=e$---as necessarily happens when $a$ is unipotent, for example---the antideveloped region $q_{{}_H}^{-1}(U)_\mathscr{G}$ has trivial holonomy.

\begin{Theorem}\label{ensnaringthm} Suppose $(\mathscr{G},\omega)$ is a Cartan geometry of type $(G,H)$ with an automorphism $\alpha\in\Aut(\mathscr{G},\omega)$ and an element $\mathscr{e}\in\mathscr{G}$ such that $\alpha(\mathscr{e})=\mathscr{e}a$ for some $a\in H$ ensnaring $U\subseteq G/H$. If $a^kga^{-k}\rightarrow e$ as $k\rightarrow+\infty$ implies $g=e$, then \smash{$\Hol_{(\hat{\zeta}_U)_\mathscr{G}(1)}\bigl(q_{{}_H}^{-1}(U)_\mathscr{G},\omega\bigr)=\{e\}$}.\end{Theorem}
\begin{proof}Suppose $\gamma \colon [0,1]\to q_{{}_H}^{-1}(U)_\mathscr{G}$ is a path such that \[\gamma(0)=\gamma(1)h^{-1}=\bigl(\hat{\zeta}_U\bigr)_\mathscr{G}(1)\] for some $h\in H$. We want to show that $\gamma_G(1)h^{-1}=e$. To do this, consider the path \[\gamma_{\zeta_U}:=(\hat{\zeta}_U)_\mathscr{G}\star\gamma\star\Rt{h}\overline{(\hat{\zeta}_U)}_\mathscr{G},\] where $\overline{\eta}$ is the \emph{reverse} of the path $\eta$, given by $\overline{\eta}(t):=\eta(1-t)$, and $\Rt{h}$ denotes right-translation by $h$.

Since \smash{$\alpha^k_*\rightarrow 0_{q_{{}_H}(\mathscr{e})}$} locally uniformly over $U_\mathscr{G}$ as $k\rightarrow +\infty$, the paths $\alpha^k q_{{}_H}(\gamma_{\zeta_U})$ converge, modulo parametrization, to $q_{{}_H}(\mathscr{e})$ to first order as before, so there exist paths $\beta_k\colon [0,1]\to H$ with $\beta_k(0)=a^k$ such that, again modulo parametrization, $\alpha^k\gamma_{\zeta_U}\beta_k^{-1}$ converges to $\mathscr{e}$ in $\Path_\mathscr{e}(\mathscr{G})$. Thus, by Lemma~\ref{shrinkingholonomy}, the developments
\[\bigl(\alpha^k\gamma_{\zeta_U}\beta_k^{-1}\bigr)_G(1) =\beta_k(0)(\gamma_{\zeta_U})_G(1)\beta_k(1)^{-1}=a^k(\gamma_{\zeta_U})_G(1)\beta_k(1)^{-1}\]
must converge to the identity element $e\in G$. Note, in particular, that the endpoint
\begin{align*}\alpha^k(\gamma_{\zeta_U}(1))\beta_k(1)^{-1}   =\alpha^k\bigl(\Rt{h}\overline{\bigl(\hat{\zeta}_U\bigr)}_\mathscr{G}(1)\bigr)\beta_k(1)^{-1}  =\alpha^k(\mathscr{e}h)\beta_k(1)^{-1}=\mathscr{e}a^kh\beta_k(1)^{-1}
\end{align*}
converges to $\mathscr{e}$, so $a^kh\beta_k(1)^{-1}\rightarrow e$, hence
\[a^k\bigl((\gamma_{\zeta_U})_G(1)h^{-1}\bigr)a^{-k} =(a^k(\gamma_{\zeta_U})_G(1)\beta_k(1)^{-1})\bigl(a^kh\beta_k(1)^{-1}\bigr)^{-1}\rightarrow ee^{-1}=e.\]
But by hypothesis, this implies that $(\gamma_{\zeta_U})_G(1)h^{-1}=e$, and
\begin{align*}(\gamma_{\zeta_U})_G(1) & =\bigl(\bigl(\hat{\zeta}_U\bigr)_\mathscr{G}\star\gamma\star\Rt{h}\overline{\bigl(\hat{\zeta}_U\bigr)}_\mathscr{G}\bigr)_G (1) =\bigl(\bigl(\hat{\zeta}_U\bigr)_\mathscr{G}\bigr)_G\!(1) \gamma_G(1)\bigl(\Rt{h}\overline{\bigl(\hat{\zeta}_U\bigr)}_\mathscr{G}\bigr)_G (1)
\\ & =\hat{\zeta}_U(1)\gamma_G(1)h^{-1}\hat{\zeta}_U(1)^{-1}h,
\end{align*} so $\gamma_G(1)h^{-1}=\hat{\zeta}_U(1)^{-1}(\gamma_{\zeta_U})_G(1)h^{-1}\hat{\zeta}_U(1)=\hat{\zeta}_U(1)^{-1}\hat{\zeta}_U(1)=e$.
\end{proof}

It is, perhaps, worth explicitly pointing out that Theorem~\ref{ensnaringthm} starkly improves previous results that exhibit local flatness near an isolated higher-order fixed point: with considerably weaker assumptions, we have still obtained a large open subset \smash{$q_{{}_H}^{-1}(U)_\mathscr{G}$} containing~$\mathscr{e}$ in its closure, and not only does the curvature vanish on this subset, but its holonomy is trivial, even without additional curvature restrictions or considerations other than the behavior of the isotropy on the model geometry. Moreover, if we also happen to know a bit about the topology of $U$, then we can even show that this subspace is isomorphic to its analogue in the model.

\begin{Corollary}\label{ensnaringthmcor} In the setting of Theorem~{\rm \ref{ensnaringthm}}, if $U$ is simply connected, then there exists a geometric embedding $\sigma\colon \bigl(q_{{}_H}^{-1}(U),\MC{G}\bigr)\hookrightarrow(\mathscr{G},\omega)$ with $\sigma\bigl(q_{{}_H}^{-1}(U)\bigr)=q_{{}_H}^{-1}(U)_\mathscr{G}$.
\end{Corollary}

\begin{proof}Naturally, we want to try to define $\sigma\colon q_{{}_H}^{-1}(U)\to q_{{}_H}^{-1}(U)_\mathscr{G}\subseteq\mathscr{G}$ by
\[\sigma(\delta(1)h):=\bigl(\hat{\zeta}_U\star\delta\bigr)_\mathscr{G}(1)h\]
for each path \smash{$\delta\in\Path_{\hat{\zeta}_U(1)}\bigl(q_{{}_H}^{-1}(U)\bigr)$} and each $h\in H$. If $\sigma$ is well defined, then it is $H$-equivariant and satisfies \smash{$\omega(\sigma_*\dot{\delta}(1))=\MC{G}(\dot{\delta}(1))$} by definition, so it is necessarily geometric. It remains, then, to prove that $\sigma$ is well defined and that~$\sigma$ is injective.

If we have two paths \smash{$\delta_1,\delta_2\in\Path_{\hat{\zeta}_U(1)}\bigl(q_{{}_H}^{-1}(U)\bigr)$} with the same endpoint $\delta_1(1)=\delta_2(1)$, then we can concatenate $\delta_1$ with the reverse of~$\delta_2$ to get a loop $\delta_1\star\bar{\delta}_2$ based at $\hat{\zeta}_U(1)$. In particular, to prove that \smash{$\bigl(\hat{\zeta}_U\star\delta_1\bigr)_\mathscr{G}(1)=\bigl(\hat{\zeta}_U\star\delta_2\bigr)_\mathscr{G}(1)$} for all such paths~$\delta_1$ and~$\delta_2$, it suffices to prove that \smash{$\bigl(\hat{\zeta}_U\star\delta\bigr)_\mathscr{G}(1)=\bigl(\hat{\zeta}_U\bigr)_\mathscr{G}(1)$} for each loop \smash{$\delta\in\Path_{\hat{\zeta}_U(1)}\bigl(q_{{}_H}^{-1}(U)\bigr)$} based at $\hat{\zeta}_U(1)$.

For each loop \smash{$\delta\in\Path_{\hat{\zeta}_U(1)}\bigl(q_{{}_H}^{-1}(U)\bigr)$}, we can find a loop $\beta\in\Path_e(H)$ such that there exists a null-homotopy \smash{$c\colon [0,1]^2\to q_{{}_H}^{-1}(U)$} from $c_0=\delta\beta^{-1}$ to the constant loop \smash{$c_1=\hat{\zeta}_U(1)$}, since $U$ is simply connected. From Lemma~\ref{snarelemma}, we get a well-defined antidevelopment \smash{$\bigl(\hat{\zeta}_U\star c_s\bigr)_\mathscr{G}$} for each $s\in[0,1]$, and we can show that the map \smash{$\sigma(c)\colon [0,1]^2\to q_{{}_H}^{-1}(U)_\mathscr{G}$} given by
\[\sigma(c)_s(t):=\sigma(c_s(t))=\bigl(\hat{\zeta}_U\star c_s\bigr)_\mathscr{G}\bigl(\tfrac{1+t}{2}\bigr)\]
satisfies $\sigma(c)^*\omega=c^*\MC{G}$, essentially for the same reason as Theorem~6.1 in \cite[Chapter~3]{Sharpe1997}: the distribution $\ker(c^*\MC{G}-\omega)$ on $[0,1]\times q_{{}_H}^{-1}(U)_\mathscr{G}$ is integrable because the curvature $\mathrm{d}\omega+\tfrac{1}{2}[\omega,\omega]$ vanishes on the subspace \smash{$q_{{}_H}^{-1}(U)_\mathscr{G}$}, so we can build the graph of $\sigma(c)$ locally from the integral submanifolds for the distribution, and \smash{$\sigma(c)_s\colon t\mapsto\bigl(\hat{\zeta}_U\star c_s\bigr)_\mathscr{G}\bigl(\tfrac{1+t}{2}\bigr)$} is the unique path starting at \smash{$\sigma(c)_s(0)=\bigl(\hat{\zeta}_U\bigr)_\mathscr{G}(1)$} such that $\sigma(c)_s^*\omega=c_s^*\MC{G}$. Since $\sigma(c)_s(0)$ and $\sigma(c)_s(1)$ are constant in $s$ and $\sigma(c)_1(t)$ is constant in $t$, it follows that $\sigma(c)_s(0)=\sigma(c)_1(t)=\sigma(c)_s(1)=\bigl(\hat{\zeta}_U\bigr)_\mathscr{G}(1)$, so $\sigma(c)$ is a null-homotopy as well. In particular,
\begin{align*}\bigl(\hat{\zeta}_U\star\delta\bigr)_\mathscr{G}(1) & =\bigl(\hat{\zeta}_U\star\delta\beta^{-1}\bigr)_\mathscr{G}(1)\beta(1)=\bigl(\hat{\zeta}_U\star\delta\beta^{-1}\bigr)_\mathscr{G}(1) \\ & =\bigl(\hat{\zeta}_U\star c_0\bigr)_\mathscr{G}(1)=\sigma(c)_0(1)=\sigma(c)_0(0)=\bigl(\hat{\zeta}_U\bigr)_\mathscr{G}(1),
\end{align*} so $\sigma$ is a well-defined geometric map.

Finally, to show that $\sigma$ is injective, we consider the shifted developing map \[\hat{\zeta}_U(1)\mathrm{dev}_{(\hat{\zeta}_U)_\mathscr{G}(1)}\colon \ q_{{}_H}^{-1}(U)_\mathscr{G}\to q_{{}_H}^{-1}(U)\] given by \smash{$\gamma(1)\mapsto\hat{\zeta}_U(1)\gamma_G(1)$} for each path \smash{$\gamma\in\Path_{(\hat{\zeta}_U)_\mathscr{G}(1)}\bigl(q_{{}_H}^{-1}(U)_\mathscr{G}\bigr)$}; by Theorem~\ref{ensnaringthm}, we know that \smash{$\Hol_{(\hat{\zeta}_U)_\mathscr{G}(1)}\bigl(q_{{}_H}^{-1}(U)_\mathscr{G},\omega\bigr)=\{e\}$}, so the map \smash{$\hat{\zeta}_U(1)\mathrm{dev}_{(\hat{\zeta}_U)_\mathscr{G}(1)}$} is necessarily well defined on \smash{$q_{{}_H}^{-1}(U)_\mathscr{G}$}. Because both $\sigma$ and \smash{$\hat{\zeta}_U(1)\mathrm{dev}_{(\hat{\zeta}_U)_\mathscr{G}(1)}$} are geometric maps, it follows that the composition \smash{$\sigma\circ\bigl(\hat{\zeta}_U(1)\mathrm{dev}_{(\hat{\zeta}_U)_\mathscr{G}(1)}\bigr)$} is a geometric map from \smash{$q_{{}_H}^{-1}(U)_\mathscr{G}$} to itself such that \[\bigl(\sigma\circ\bigl(\hat{\zeta}_U(1)\mathrm{dev}_{(\hat{\zeta}_U)_\mathscr{G}(1)}\bigr)\bigr) \bigl(\bigl(\hat{\zeta}_U\bigr)_\mathscr{G}(1)\bigr)=\sigma\bigl(\hat{\zeta}_U(1)\bigr) =\bigl(\hat{\zeta}_U\bigr)_\mathscr{G}(1),\] so \smash{$\hat{\zeta}_U(1)\mathrm{dev}_{(\hat{\zeta}_U)_\mathscr{G}(1)}=\sigma|_{q_H^{-1}(U)_\mathscr{G}}^{-1}$}, hence $\sigma$ is injective. \end{proof}

In the next section, we will further refine Corollary~\ref{ensnaringthmcor} in the case where the open set $U=(G/H)\setminus\mathrm{Fix}_{G/H}(a)$ is simply connected and $\mathrm{Fix}_{G/H}(a)$ has codimension at least two. For such situations, we will see that \smash{$q_{{}_H}^{-1}(U)_\mathscr{G}$} already makes up so much of $\mathscr{G}$ that we can essentially just extend~$\sigma^{-1}$ to all of~$\mathscr{G}$ by filling in the small gaps where the fixed points would go.

\section{Applications}\label{applications}
Here, we will prove Theorems \ref{thmA}, \ref{thmB}, and \ref{thmC} from the introduction. To do this, we will prove a more general result, and then we will show that the hypotheses of this result are satisfied in each of the desired cases.

Throughout this section, we will need a notion of codimension; it will be convenient to use the following terminology, which comes from considerations of transversality.

\begin{Definition}For a smooth manifold $M$ and a closed subset $N\subseteq M$, we will say that $N$ has \emph{codimension at least two} if and only if, for each path $\gamma \colon [0,1]\to M$ with $\gamma(0)\not\in N$, there exists a homotopy $c\colon [0,1]^2\to M$, $(s,t)\mapsto c_s(t)$ such that $c_0=\gamma$, $c_s(0)=\gamma(0)$ for all $s\in[0,1]$, and $c((0,1]\times[0,1])\cap N=\varnothing$.\end{Definition}

Closed submanifolds and subvarieties $N\subseteq M$ will, of course, have codimension at least two in this sense if $\dim(M)-\dim(N)\geq 2$. Moreover, the complement of~$N$ must be dense in~$M$ if~$N$ satisfies this definition: for each point $p\in M$, we can find a sequence of points of the form $c_{s_i}(t_i)\not\in N$ with $s_i\rightarrow 0$ such that $c_{s_i}(t_i)\rightarrow p$.

\subsection{General result}
In this subsection, we will prove the following result, from which Theorems \ref{thmA}, \ref{thmB}, and \ref{thmC} will follow.

\begin{Theorem}\label{ensnaringthm+} Let $(\mathscr{G},\omega)$ be a Cartan geometry of type $(G,H)$ over a connected smooth mani\-fold~$M$, $\alpha\in\Aut(\mathscr{G},\omega)$, and $\mathscr{e}\in\mathscr{G}$ such that $\alpha(\mathscr{e})=\mathscr{e}a$ for some $a\in H$. Suppose that the isotropy $a\in H$ ensnares $U=(G/H)\setminus\mathrm{Fix}_{G/H}(a)$, that $U$ is simply connected, and that $\mathrm{Fix}_{G/H}(a)$ has codimension at least two in $G/H$. If $a^kga^{-k}\rightarrow e$ as $k\rightarrow +\infty$ implies $g=e$, then $(\mathscr{G},\omega)$ geometrically embeds onto a dense open subset of the Klein geometry $(G,\MC{G})$ of type $(G,H)$.\end{Theorem}

By Corollary~\ref{ensnaringthmcor}, the hypotheses of the theorem tell us that there is a geometric embedding $\sigma\colon \bigl(q_{{}_H}^{-1}(U),\MC{G}\bigr)\hookrightarrow (\mathscr{G},\omega)$. It turns out that this is already sufficient to prove the theorem: $U$ has complement of codimension at least two, so we can use \cite[Theorem~1.8]{Frances2008} to get a geometric isomorphism between $(\mathscr{G},\omega)$ and a dense open subset of $G$. For completeness---and to continue to demonstrate the utility of holonomy for obtaining global results on Cartan geometries---we provide an alternative proof.

Intuitively, the image $\sigma(U)=U_\mathscr{G}$ of $\sigma$ should already make up nearly all of $M$, and the holonomy is trivial over this subset. The main idea behind our proof of Theorem~\ref{ensnaringthm+}, then, is to use the codimension of the complement of $U$ in $G/H$ to homotope loops in $M$ an arbitrarily small amount to get loops in $\sigma(U)$, thus forcing the holonomy to be globally trivial.

\begin{figure}[t]
\centering
\includegraphics[width=0.4\textwidth]{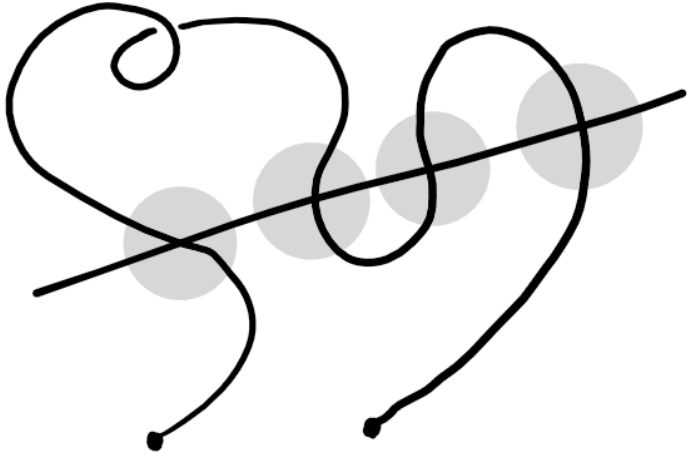}
\put(-122,117){\makebox(0,0)[lb]{\small $\gamma_G=c_0$}}
\put(2,90){\makebox(0,0)[lb]{\small $G\setminus\mathscr{F}$}}

\caption{The path $\gamma_G$ in $G$, with neighborhoods of its intersection with $G\setminus\mathscr{F}$ highlighted in gray.}\label{codim2fig0}
\end{figure}

\begin{figure}[t]
\centering
\includegraphics[width=0.4\textwidth]{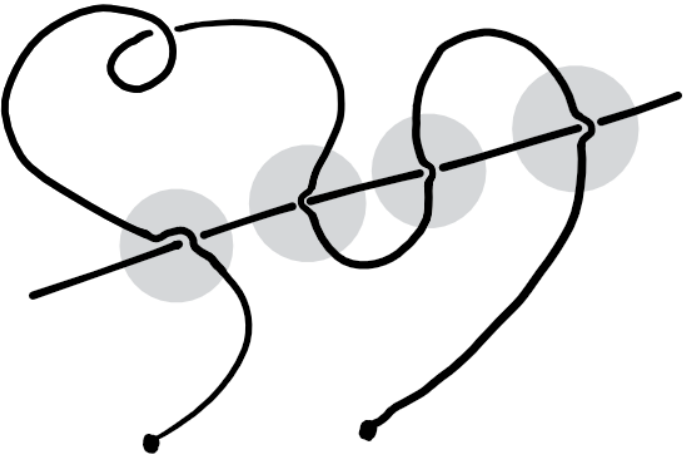}
\put(-122,118){\makebox(0,0)[lb]{\small $c_\varepsilon$}}
\put(2,90){\makebox(0,0)[lb]{\small $G\setminus\mathscr{F}$}}

\caption{A path $c_\varepsilon$ in $G$, which satisfies $c_\varepsilon([0,1])\subseteq\mathscr{F}$ and coincides with $\gamma_G=c_0$ outside of the highlighted neighborhoods of the intersections of $\gamma_G$ with $G\setminus\mathscr{F}$.}
\label{codim2fig1}
\end{figure}

\begin{Proposition}\label{codim2}
Suppose $(\mathscr{G},\omega)$ is a Cartan geometry of type $(G,H)$ over a connected smooth manifold $M$. Let $\mathscr{F}\subseteq G$ be an $H$-invariant open subset of $G$ such that the complement $G\setminus\mathscr{F}$ has codimension at least two. If there exists a geometric embedding $\sigma\colon (\mathscr{F},\MC{G})\hookrightarrow (\mathscr{G},\omega)$, then $\Hol_\mathscr{e}(\mathscr{G},\omega)=\{e\}$ for any $\mathscr{e}\in\mathscr{G}$.
\end{Proposition}
\begin{proof}Suppose $\gamma \colon [0,1]\to\mathscr{G}$ is the lift of a loop $q_{{}_H}(\gamma)$ on $M$ with $\gamma(0)=\mathscr{e}$, and $h_\gamma\in H$ is such that $\gamma(0)h_\gamma=\gamma(1)$. We want to prove that $\gamma_G(1)h_\gamma^{-1}=e$; since $M$ is connected, we may assume without loss of generality that $\mathscr{e}\in\sigma(\mathscr{F})$, and since $G$ acts transitively on itself, we may further assume that $e\in\mathscr{F}$ and $\sigma(e)=\mathscr{e}$.

Because $G\setminus\mathscr{F}$ has codimension at least two in $G$, there exists a homotopy $c\colon [0,1]^2\to G$ such that $c_0=\gamma_G$, $c_s(0)=\gamma_G(0)=e$ for all $s\in[0,1]$, and $c((0,1]\times[0,1])\subseteq\mathscr{F}$; see Figures~\ref{codim2fig0} and~\ref{codim2fig1} for illustrations. Then, since $\sigma$ is a geometric map and $c_s(0)=e\in\mathscr{F}$ for each $s\in[0,1]$, $\sigma(c_s)_G=c_s(0)^{-1}c_s=c_s$ for all $s\in (0,1]$, hence---by Lemma~\ref{cdop}---$\sigma(c_s)$ converges to $\gamma\in\Path_\mathscr{e}(\mathscr{G})$ as $s\rightarrow 0$. In particular, \[\sigma\bigl(c_s(1)h_\gamma^{-1}\bigr)=\sigma(c_s(1))h_\gamma^{-1}\rightarrow\gamma(1)h_\gamma^{-1}=\gamma(0)=\sigma(e),\] so since $\sigma$ is an embedding, we must have $c_s(1)h_\gamma^{-1}\rightarrow\gamma_G(1)h_\gamma^{-1}=e$, hence $\Hol_\mathscr{e}(\mathscr{G},\omega)=\{e\}$.
\end{proof}

Recall that, whenever $\Hol_\mathscr{e}(\mathscr{G},\omega)=\{e\}$, the developing map \[\mathrm{dev}_\mathscr{e}\colon \ (\mathscr{G},\omega)\to(G,\MC{G})\] based at $\mathscr{e}$ is a~well-defined geometric map given by $\gamma(1)h\mapsto\gamma_G(1)h$ for every $h\in H$ and piecewise smooth path $\gamma \colon [0,1]\to\mathscr{G}$ starting at $\gamma(0)=\mathscr{e}$. Therefore, if we can show that $\mathrm{dev}_\mathscr{e}$ is injective, then the proof of Theorem~\ref{ensnaringthm+} will be done.

To do this, we will use the following lemma, which tells us that a geometric map (or more generally any local homeomorphism from a Hausdorff space) is injective when its restriction to the preimage of a dense subset is injective. We imagine that this lemma is well known, but could not find a reference.

\begin{Lemma}\label{injpreimage} Suppose $\psi\colon (\mathscr{G},\omega)\to(\mathscr{Q},\upsilon)$ is a geometric map and that $V\subseteq\mathscr{Q}$ is a dense subset. If $\psi|_{\psi^{-1}(V)}$ is injective, then $\psi$ is injective.
\end{Lemma}
\begin{proof}Suppose that $\psi(\mathscr{g})=\psi(\mathscr{g}')$ for some $\mathscr{g},\mathscr{g}'\in\mathscr{G}$. Since $\psi$ is a local diffeomorphism, there exists an open neighborhood $W$ of $\mathscr{g}'$ such that $\psi|_W$ is a diffeomorphism onto its image. If $(\mathscr{g}_i)$ is a sequence in $\psi^{-1}(V)$ converging to $\mathscr{g}$, which must exist because density of $V$ in $\mathscr{Q}$ implies density of $\psi^{-1}(V)$ in $\mathscr{G}$, then $\psi(\mathscr{g}_i)\rightarrow\psi(\mathscr{g})=\psi(\mathscr{g}')$, so $\psi|_W^{-1}(\psi(\mathscr{g}_i))\rightarrow\psi|_W^{-1}(\psi(\mathscr{g}'))=\mathscr{g}'$. But $\psi|_{\psi^{-1}(V)}$ is injective, so $\psi|_W^{-1}(\psi(\mathscr{g}_i))=\mathscr{g}_i\rightarrow\mathscr{g}$; since $\mathscr{G}$ is Hausdorff, this means $\mathscr{g}=\mathscr{g}'$.
\end{proof}

\begin{Corollary}\label{codim2cor} In the setting of Proposition~{\rm \ref{codim2}}, the developing map $\mathrm{dev}_\mathscr{e}\colon (\mathscr{G},\omega)\to(G,\MC{G})$ is a geometric embedding.
\end{Corollary}
\begin{proof}As we did in the proof of Proposition~\ref{codim2}, we may assume that $e\in\mathscr{F}$ and $\sigma(e)=\mathscr{e}$. By definition, we have that $\mathrm{dev}_\mathscr{e}(\gamma)=\gamma_G$ whenever $\gamma(0)=\mathscr{e}$, so $(\mathrm{dev}_\mathscr{e}\circ\sigma)(\gamma_G)=\mathrm{dev}_\mathscr{e}(\gamma)=\gamma_G$ whenever $\gamma_G$ is contained in $\mathscr{F}$. Hence, $\mathrm{dev}_\mathscr{e}\circ\sigma$ is the identity map on $\mathscr{F}$.

It just remains to show that $\mathrm{dev}_\mathscr{e}\colon \mathscr{G}\to G$ is injective. To this end, suppose $\mathrm{dev}_\mathscr{e}(\mathscr{g})\in\mathscr{F}$, and let $\gamma\in\Path_\mathscr{e}(\mathscr{G})$ be a path such that $\gamma(1)\in\mathscr{g}H$. Then, $\mathrm{dev}_\mathscr{e}(\gamma)$ is a path in $G$ with $\mathrm{dev}_\mathscr{e}(\gamma(0))=\mathrm{dev}_\mathscr{e}(\mathscr{e})=e$ and $\mathrm{dev}_\mathscr{e}(\gamma(1))=\gamma_G(1)\in\mathrm{dev}_\mathscr{e}(\mathscr{g})H$. Since $G\setminus\mathscr{F}$ has codimension at least two, there exists a homotopy $c\colon [0,1]^2\to G$ such that $c_0=\mathrm{dev}_\mathscr{e}(\gamma)=\gamma_G$, $c_s(0)=e$ for all $s\in[0,1]$, and $c((0,1]\times[0,1])\subseteq\mathscr{F}$. By continuity, $\sigma(c_s(1))\rightarrow \sigma(c_0(1))=\sigma(\gamma_G(1))$ as $s\rightarrow 0$, and by Lemma~\ref{cdop}, $\sigma(c_s(1))\rightarrow \gamma(1)$ as $s\rightarrow 0$. Because $\mathscr{G}$ is Hausdorff, this means that $\gamma(1)=\sigma(\gamma_G(1))\in\sigma(\mathscr{F})$, hence that $\mathscr{g}\in\sigma(\mathscr{F})$. Thus, $\mathrm{dev}_\mathscr{e}^{-1}(\mathscr{F})=\sigma(\mathscr{F})$, so because $\mathrm{dev}_\mathscr{e}|_{\sigma(\mathscr{F})}$ is injective, $\mathrm{dev}_\mathscr{e}$ is injective by Lemma~\ref{injpreimage}.
\end{proof}

In summary, under the hypotheses of Theorem~\ref{ensnaringthm+}, we have shown that the developing map $\mathrm{dev}_\mathscr{e}$ gives a geometric embedding from $(\mathscr{G},\omega)$ onto a dense open subset of the Klein geometry $(G,\MC{G})$, so the proof is complete.

\subsection{Proof of Theorem~\ref{thmA}} Recall that, in this case, our model geometry is $(\PGL_{m+1}\mathbb{C},P)$, which corresponds to complex projective geometry over $\mathbb{CP}^m$. We want to prove that, given a Cartan geometry $(\mathscr{G},\omega)$ of type $(\PGL_{m+1}\mathbb{C},P)$ over a connected smooth manifold $M$ and a nontrivial geometric automorphism $\alpha\in\Aut(\mathscr{G},\omega)$ such that $\alpha(\mathscr{e})\in\mathscr{e}P_+$ for some $\mathscr{e}\in\mathscr{G}$, there exists a geometric embedding of $(\mathscr{G},\omega)$ onto a dense open subset of the Klein geometry $(\PGL_{m+1}\mathbb{C},\MC{\PGL_{m+1}\mathbb{C}})$.

For $\alpha(\mathscr{e})=\mathscr{e}a$, let us write \[a=\begin{pmatrix}1 & \beta \\ 0 & \mathds{1}\end{pmatrix}\] for some $\beta^\top\in\mathbb{C}^m$. In view of Theorem~\ref{ensnaringthm+}, it suffices to prove that $a$ ensnares the complement $U=\mathbb{CP}^m\setminus\mathrm{Fix}_{\mathbb{CP}^m}(a)$, that $U$ is simply connected, and that $\mathrm{Fix}_{\mathbb{CP}^m}(a)$ has codimension at least two.

For arbitrary \[q_{{}_P}(g)=\begin{pmatrix}r \\ x\end{pmatrix}\in\mathbb{CP}^m,\] with $r\in\mathbb{C}$ and $x\in\mathbb{C}^m$, we have \[a^k(q_{{}_P}(g))=a^k\begin{pmatrix}r \\ x\end{pmatrix}=\begin{pmatrix}1 & k\beta \\ 0 & \mathds{1}\end{pmatrix}\begin{pmatrix}r \\ x\end{pmatrix}=\begin{pmatrix}r+k\beta(x) \\ x\end{pmatrix},\] which is equal to \[\begin{pmatrix}1 \\ \dfrac{1}{r+k\beta(x)}x\end{pmatrix}\] whenever $r+k\beta(x)\neq 0$. In particular, fixed points are precisely those points for which $\beta(x)=0$, and if $\beta(x)\neq 0$, then $a^k(q_{{}_P}(g))$ must go to $q_{{}_P}(e)$ as $k\rightarrow +\infty$. Since
\[\mathrm{Fix}_{\mathbb{CP}^m}(a)=\left\{\begin{pmatrix}r \\ x\end{pmatrix}\in\mathbb{CP}^m\mid \beta(x)=0\right\}\]
is a complex projective hyperplane in $\mathbb{CP}^m$, $U=\mathbb{CP}^m\setminus\mathrm{Fix}_{\mathbb{CP}^m}(a)$ is diffeomorphic to $\mathbb{C}^m$, hence it is simply connected. Moreover, since $\beta(x)=0$ implies $\mathrm{Re}(\beta(x))=0$ and $\mathrm{Im}(\beta(x))=0$, $\mathrm{Fix}_{\mathbb{CP}^m}(a)$ has (real) codimension 2 in $\mathbb{CP}^m$.

Now, all that is left to do is prove that $a$ ensnares $U$. For this, note that $\zeta_U\colon [0,1]\to\mathbb{CP}^m$, given by \smash{$t\mapsto\bigl(\begin{smallmatrix}1 \\ ty\end{smallmatrix}\bigr)$} for some $y\in\mathbb{C}^m$ such that $\beta(y)\in(0,+\infty)\subset\mathbb{R}$, shrinks into itself under forward iterates of $a$, so that $a^k(\zeta_U(t))=\smash{\zeta_U\bigl(\tfrac{t}{1+k\beta(y)t}\bigr)}$. Thus, for \smash{$\tau_k\colon t\mapsto\tfrac{t}{1+k\beta(y)(1-t)}$}, we get $a^k(\zeta_U(\tau_k(t)))=\smash{\zeta_U\bigl(\tfrac{t}{1+k\beta(y)}\bigr)}$, hence $a^k(\zeta_U\circ\tau_k)$ converges to $q_{{}_P}(e)$ to first order.
Moreover, since $a^k(q_{{}_P}(g))$ converges to $q_{{}_P}(e)$ for each $q_{{}_P}(g)\in U$, we may assume that $r+k\beta(x)\neq 0$ for all $k\geq 0$ after replacing $q_{{}_P}(g)$ with $a^n(q_{{}_P}(g))$ for suitably large $n$. This allows us to work in the coordinate chart $r=1$, and we know that in those coordinates, $a^k$ is given by \smash{$x\mapsto\tfrac{1}{1+k\beta(x)}x$}. Therefore, for $\xi\in T_{q_P(g)}U$ and considering $a^k_*\xi$ as a point derivation on the vector of coordinate functions~$x$,
\begin{align*}
\bigl(a^k_*\xi\bigr)(x) & =\xi\left(\frac{1}{1+k\beta(x)}x\right)=\frac{1}{1+k\beta(x)}\left(\xi(x)-\frac{k\xi(\beta(x))}{1+k\beta(x)}x\right) \\ & =\frac{1}{(1+k\beta(x))^2}\big(\xi(x)+k(\beta(x)\xi(x)-\xi(\beta(x))x)\big),
\end{align*}
so since $(1+k\beta(x))^2$ grows quadratically in $k$ while the numerator $\xi(x)+k(\beta(x)\xi(x)-\xi(\beta(x))x)$ only grows linearly in $k$, $a^k_*\xi\rightarrow 0_{q_P(e)}$ as $k\rightarrow+\infty$. Indeed, for all sufficiently small open subsets $V\subset TU\subseteq T(G/P)$, $\pi\colon T(G/P)\to G/P$ the natural projection from the tangent bundle to $G/P\cong\mathbb{CP}^m$, $\|\cdot\|$ a norm on $\mathbb{C}^m$, and $k>(\inf_{x\in\pi(V)}|\beta(x)|)^{-1}$, we have
\[
\left\|\frac{1}{1+k\beta(x)}x\right\|=\frac{\|x\|}{|1+k\beta(x)|} \leq\frac{\sup_{x\in\pi(V)}\|x\|}{\inf_{x\in\pi(V)}|\beta(x)|k-1}\rightarrow 0\]
and
\begin{align*}
\|(a^k_*\xi)(x)\| ={}&\left\|\xi\left(\frac{1}{1+k\beta(x)}x\right)\right\| \\ \leq{} &\frac{1}{(\inf_{x\in\pi(V)}|\beta(x)|k-1)^2}\Bigl(k\bigl(\sup_{x\in\pi(V)}|\beta(x)|\sup_{\xi\in V}\|\xi(x)\| \\ & +\sup_{\xi\in V}|\xi(\beta(x))|\sup_{x\in\pi(V)}\|x\|\bigr)+\sup_{\xi\in V}\|\xi(x)\|\Bigr)\rightarrow 0
\end{align*}
for $x$ in $\pi(V)$ and $\xi\in V$ over the chart $r=1$, so the convergence of $a^k_*$ to $0_{q_{{}_P}(e)}$ is locally uniform over $U$. In other words, $a$ ensnares~$U$, so the proof of Theorem~\ref{thmA} is complete.

\subsection{Proof of Theorem~\ref{thmB}}
This proof is largely the same as the one for Theorem~\ref{thmA}. In this case, the model becomes $(\PGL_{m+1}\mathbb{H},P)$, with $\PGL_{m+1}\mathbb{H}/P\cong\mathbb{HP}^m$. We want to prove, given a Cartan geometry $(\mathscr{G},\omega)$ of type $(\PGL_{m+1}\mathbb{H},P)$ over a connected smooth manifold $M$ and a nontrivial automorphism $\alpha\in\Aut(\mathscr{G},\omega)$ such that $\alpha(\mathscr{e})\in\mathscr{e}P_+$ for some $\mathscr{e}\in\mathscr{G}$, that there exists a geometric embedding of $(\mathscr{G},\omega)$ onto a dense open subset of the Klein geometry $(\PGL_{m+1}\mathbb{H},\MC{\PGL_{m+1}\mathbb{H}})$.

As before, for $\alpha(\mathscr{e})=\mathscr{e}a$, we write \[a=\begin{pmatrix}1 & \beta \\ 0 & \mathds{1}\end{pmatrix}\] for some $\beta^\top\in\mathbb{H}^m$. Again, our goal is just to apply Theorem~\ref{ensnaringthm+}, so we just need to prove that $a$ ensnares $U=\mathbb{HP}^m\setminus\mathrm{Fix}_{\mathbb{HP}^m}(a)$, that $U$ is simply connected, and that $\mathrm{Fix}_{\mathbb{HP}^m}(a)$ has codimension at least two.

For arbitrary \[q_{{}_P}(g)=\begin{pmatrix}r \\ x\end{pmatrix}\in\mathbb{HP}^m,\] with $r\in\mathbb{H}$ and $x\in\mathbb{H}^m$, we have \[a^k(q_{{}_P}(g))=a^k\begin{pmatrix}r \\ x\end{pmatrix}=\begin{pmatrix}1 & k\beta \\ 0 & \mathds{1}\end{pmatrix}\begin{pmatrix}r \\ x\end{pmatrix}=\begin{pmatrix}r+k\beta(x) \\ x\end{pmatrix},\] which is equal to \[\begin{pmatrix}1 \\ x(r+k\beta(x))^{-1}\end{pmatrix}\] whenever $r+k\beta(x)\neq 0$. In particular, fixed points are precisely those points for which $\beta(x)=0$, and if $\beta(x)\neq 0$, then $a^k(q_{{}_P}(g))$ must go to $q_{{}_P}(e)$ as $k\rightarrow +\infty$. This time, the set $\mathrm{Fix}_{\mathbb{HP}^m}(a)$ of all points for which $\beta(x)=0$ is of codimension $\dim(\mathbb{H})=4$, since $\beta(x)=0$ implies its projections to the real subspaces of $\mathbb{H}$ generated by $1$, $\mathrm{i}$, $\mathrm{j}$, and $\mathrm{i}\mathrm{j}$ all vanish as well. Since $\mathbb{HP}^m$ is simply connected and $\mathrm{Fix}_{\mathbb{HP}^m}(a)$ has codimension greater than two, it also follows that $U=\mathbb{HP}^m\setminus\mathrm{Fix}_{\mathbb{HP}^m}(a)$ is simply connected.

We again end the proof by proving that $a$ ensnares $U$. As before, the path $\zeta_U\colon [0,1]\to\mathbb{HP}^m$ given by $t\mapsto\bigl(\begin{smallmatrix}1 \\ ty\end{smallmatrix}\bigr)$ for some $y\in\mathbb{H}^m$ such that $\beta(y)\in(0,+\infty)\subset\mathbb{R}$, shrinks into itself under forward iterates of $a$, so $a^k(\zeta_U)$ must converge to $q_{{}_P}(e)$ to first order modulo parametrization. After replacing $q_{{}_P}(g)\in U$ with some forward iterate under $a$, we can assume $r+k\beta(x)\neq 0$ for all $k\geq 0$, allowing us to work in the coordinate chart $r=1$, where we know $a^k$ is given by $x\mapsto x(1+k\beta(x))^{-1}$. Considering $a^k_*\xi$ as a point derivation on the coordinate vector $x$,
\begin{align*}\bigl(a^k_*\xi\bigr)(x)   =\xi\bigl(x(1+k\beta(x))^{-1}\bigr)  =\big(\xi(x)-kx(1+k\beta(x))^{-1}\xi(\beta(x))\big)(1+k\beta(x))^{-1},
\end{align*} which goes to $0$ as $k\rightarrow+\infty$, so $a^k_*(\xi)\rightarrow 0_{q_{{}_P}(e)}$. Indeed, since the convergence comes from the disparity of degrees in the polynomials in $k$ as before, we again get that $a^k_*$ converges locally uniformly to $0_{q_{{}_P}(e)}$, hence $a$ ensnares $U$ and Theorem~\ref{thmB} is proven.

\subsection{Proof of Theorem~\ref{thmC}} This time, our model geometry is given by $(\PU(\mathrm{h}_{p,q}),P)$, basically corresponding to the natural CR structure on the null-cone for $\mathrm{h}_{p,q}$ in $\mathbb{CP}^{p+q+1}$. Our goal is to prove that, if $(\mathscr{G},\omega)$ is a Cartan geometry of type $(\PU(\mathrm{h}_{p,q}),P)$ over a connected smooth manifold $M$ and $\alpha\in\Aut(\mathscr{G},\omega)$ is an automorphism such that $\alpha(\mathscr{e})=\mathscr{e}a$ for some non-null $a\in P_+$, then there exists a geometric embedding from $(\mathscr{G},\omega)$ onto a dense open subset of the Klein geometry $(\PU(\mathrm{h}_{p,q}),\MC{\PU(\mathrm{h}_{p,q})})$. Again, this comes down to applying Theorem~\ref{ensnaringthm+}.

Let us write \[a=\begin{pmatrix}1 & \beta & \mathrm{i}s-\tfrac{1}{2}\beta I_{p,q} \bar{\beta}^\top \vspace{1mm} \\ 0 & \mathds{1} & -I_{p,q}\bar{\beta}^\top \\ 0 & 0 & 1\end{pmatrix},\] where, by hypothesis, $\beta I_{p,q}\bar{\beta}^\top\neq 0$ if $\beta\neq 0$. For arbitrary \[q_{{}_P}(g)=\begin{pmatrix}r \\ x \\ c\end{pmatrix}\in\mathrm{Null}(\mathrm{h}_{p,q})\subseteq\mathbb{CP}^{p+q+1},\] with $r,c\in\mathbb{C}$ and $x\in\mathbb{C}^{p+q}$ such that $2\mathrm{Re}(\bar{r}c)+\bar{x}^\top I_{p,q}x=0$, we get
\begin{align*}a^k(q_{{}_P}(g)) & =a^k\begin{pmatrix}r \\ x \\ c\end{pmatrix}=\begin{pmatrix}1 & k\beta & \mathrm{i}ks-\tfrac{k^2}{2}\beta I_{p,q} \bar{\beta}^\top \vspace{1mm}\\ 0 & \mathds{1} & -kI_{p,q}\bar{\beta}^\top \\ 0 & 0 & 1\end{pmatrix}\begin{pmatrix}r \\ x \\ c\end{pmatrix} \\ & =\begin{pmatrix}r+k(\beta(x)+\mathrm{i}cs)-\tfrac{ck^2}{2}\beta I_{p,q}\bar{\beta}^\top \vspace{1mm}\\ x-kcI_{p,q}\bar{\beta}^\top \\ c\end{pmatrix},
\end{align*} so $q_{{}_P}(g)\in\mathrm{Fix}_{\mathrm{Null}(\mathrm{h}_{p,q})}(a)$ if and only if $\beta(x)+\mathrm{i}cs=0$ and $cI_{p,q}\bar{\beta}^\top=0$.
If $\beta=0$ and $s\neq 0$, then this is just the set of points such that $c=0$, which has codimension $\dim(\mathbb{C})=2$, and $U=\mathrm{Null}(\mathrm{h}_{p,q})\setminus\mathrm{Fix}_{\mathrm{Null}(\mathrm{h}_{p,q})}(a)$ is diffeomorphic to $\mathrm{i}\mathbb{R}\times\mathbb{C}^{p+q}$ by normalizing to $c=1$. If $\beta\neq 0$, then we get the set of points such that $\beta(x)=0$ and $c=0$, which has codimension 4. Either way, $\mathrm{Fix}_{\mathrm{Null}(\mathrm{h}_{p,q})}(a)$ has codimension at least two, $U=\mathrm{Null}(\mathrm{h}_{p,q})\setminus\mathrm{Fix}_{\mathrm{Null}(\mathrm{h}_{p,q})}(a)$ is simply connected, and $a^k(q_{{}_P}(g))$ converges to $q_{{}_P}(e)$ as $k\rightarrow+\infty$ whenever $q_{{}_P}(g)$ is not a fixed point.

It remains to prove that $a$ ensnares $U$. In this case, we can choose $\zeta_U\colon [0,1]\to\mathrm{Null}(\mathrm{h}_{p,q})$ to be given by
\[\zeta_U(t):=\begin{pmatrix}t(t-s)-\tfrac{\mathrm{i}}{2}\beta I_{p,q}\bar{\beta}^\top \vspace{1mm}\\ -t\mathrm{i}I_{p,q}\bar{\beta}^\top \vspace{1mm}\\ t^2\mathrm{i}\end{pmatrix},\]
so that \smash{$a^k(\zeta_U(t))=\zeta_U\bigl(\tfrac{t}{1+kt}\bigr)$} and hence $a^k(\zeta_U(\tau_k(t)))=\zeta_U(\tfrac{t}{1+k})$ for the reparametrizations \smash{$\tau_k\colon t\mapsto\tfrac{t}{1+k(1-t)}$} of $[0,1]$. From this, it follows that the sequence $a^k(\zeta_U\circ\tau_k)$ converges to $q_{{}_P}(e)$ to first order. After possibly replacing $q_{{}_P}(g)$ with some forward iterate under $a$, we can assume that \smash{$r+k(\beta(x)+\mathrm{i}cs)-\tfrac{ck^2}{2}\beta I_{p,q}\bar{\beta}^\top\neq 0$} for all $k\geq 0$, so we can work in the chart where we normalize to $r=1$, with coordinates given by \[\begin{pmatrix}1 \\ x \\ c\end{pmatrix}=\begin{pmatrix}1 \\ x \\ \mathrm{i}\,\mathrm{Im}(c)-\tfrac{\bar{x}^\top I_{p,q}x}{2}\end{pmatrix}\mapsto\begin{bmatrix}x \\ \mathrm{Im}(c)\end{bmatrix}\in\mathbb{C}^{p+q}\times\mathbb{R}.\]
For $r_k:=1+k(\beta(x)+\mathrm{i}cs)-\tfrac{ck^2}{2}\beta I_{p,q}\bar{\beta}^\top$, the expression for $a^k$ in these coordinates is \[\begin{bmatrix}x \\ \mathrm{Im}(c)\end{bmatrix}\mapsto\begin{bmatrix}\tfrac{1}{r_k}\bigl(x-kcI_{p,q}\bar{\beta}^\top\bigr) \\ \mathrm{Im}\bigl(\tfrac{c}{r_k}\bigr)\end{bmatrix}.\] Therefore, since
\[\xi\left(\frac{1}{r_k}\bigl(x-kcI_{p,q}\bar{\beta}^\top\bigr)\right)
=\frac{1}{r_k^2}\bigl(r_k\bigl(\xi(x)-k\xi(c)I_{p,q}\bar{\beta}^\top\bigr) -\xi(r_k)\bigl(x-kcI_{p,q}\bar{\beta}^\top\bigr)\bigr)\] goes to $0\in\mathbb{C}^{p+q}$ and \smash{$\xi\bigl(\tfrac{c}{r_k}\bigr)=\tfrac{\xi(c)r_k-c\xi(r_k)}{r_k^2}$} goes to $0\in\mathbb{C}$ as $k\rightarrow+\infty$ for each $\xi\in T_{q_P(g)}U$, we have $a^k_*(\xi)\rightarrow 0_{q_{{}_P}(e)}$. Again, the convergence to $0_{q_{{}_P}(e)}$ comes from a disparity in the degrees of the polynomials in $k$, and it follows that $a^k_*$ converges locally uniformly over~$U$, hence $a$ ensnares~$U$, which completes the proof of Theorem~\ref{thmC}.

\subsection*{Acknowledgements}
We would specifically like to thank Karin Melnick for several helpful suggestions and conversations on this topic; she was arguably instrumental in the beginning stages of this work. We are also deeply grateful to Charles Frances, both for inspiring this paper and for tolerating and tempering the author's enthusiasm, and an anonymous referee who reviewed an earlier version of this paper and managed to provide a number of useful improvements. The reviewers at SIGMA also suggested some implemented changes, and we appreciate them for their efforts. Finally, we gratefully acknowledge partial support from the Brin Graduate Fellowship at the University of Maryland and the National Science Foundation under Grant No.\ 2203493 during the time in which this research was conducted.

\pdfbookmark[1]{References}{ref}
\LastPageEnding


\begin{thebibliography}{99}
\footnotesize\itemsep=0pt

\bibitem{cprojCEMN}
Calderbank D.M.J., Eastwood M.G., Matveev V.S., Neusser K., C-projective
 geometry, \href{https://doi.org/10.1090/memo/1299}{\textit{Mem. Amer. Math.
 Soc.}} \textbf{267} (2020), v+137~pages,
 \href{https://arxiv.org/abs/1512.04516}{arXiv:1512.04516}.

\bibitem{CapMelnick2013}
{\v{C}}ap A., Melnick K., Essential {K}illing fields of parabolic geometries,
 \href{https://doi.org/10.1512/iumj.2013.62.5166}{\textit{Indiana Univ.
 Math.~J.}} \textbf{62} (2013), 1917--1953,
 \href{https://arxiv.org/abs/1208.5510}{arXiv:1208.5510}.

\bibitem{CapSlovakPG1}
{\v{C}}ap A., Slov\'ak J., Parabolic geometries.~{I}. Background and general
 theory, \textit{Math. Surveys Monogr.}, Vol.~154,
 \href{https://doi.org/10.1090/surv/154}{American Mathematical Society},
 Providence, RI, 2009.

\bibitem{HCartan1971}
Cartan H., Differential calculus, Hermann, Paris, 1971.

\bibitem{HolonomyPaper}
Erickson J.W., Intrinsic holonomy and curved cosets of {C}artan geometries,
 \href{https://doi.org/10.1007/s40879-022-00535-7}{\textit{Eur.~J. Math.}}
 \textbf{8} (2022), 446--474,
 \href{https://arxiv.org/abs/2109.05350}{arXiv:2109.05350}.

\bibitem{SprawlPaper}
Erickson J.W., A method for determining {C}artan geometries from the local
 behavior of automorphisms,
 \href{https://arxiv.org/abs/2303.00561}{arXiv:2303.00561}.

\bibitem{Frances2007}
Frances C., Sur le groupe d'automorphismes des g\'eom\'etries paraboliques de
 rang 1, \href{https://doi.org/10.1016/j.ansens.2007.07.003}{\textit{Ann. Sci.
 \'Ecole Norm. Sup.~(4)}} \textbf{40} (2007), 741--764.

\bibitem{Frances2008}
Frances C., Rigidity at the boundary for conformal structures and other
 {C}artan geometries, \href{https://arxiv.org/abs/0806.1008}{arXiv:0806.1008}.

\bibitem{Frances2012}
Frances C., Local dynamics of conformal vector fields,
 \href{https://doi.org/10.1007/s10711-011-9619-7}{\textit{Geom. Dedicata}}
 \textbf{158} (2012), 35--59,
 \href{https://arxiv.org/abs/0909.0044}{arXiv:0909.0044}.

\bibitem{Goldman2022}
Goldman W.M., Geometric structures on manifolds, \textit{Grad. Stud. Math.},
 Vol.~227, \href{https://doi.org/10.1090/gsm/227}{American Mathematical
 Society}, Providence, RI, 2022.

\bibitem{KruglikovThe2018}
Kruglikov B., The D., Jet-determination of symmetries of parabolic geometries,
 \href{https://doi.org/10.1007/s00208-017-1545-z}{\textit{Math. Ann.}}
 \textbf{371} (2018), 1575--1613,
 \href{https://arxiv.org/abs/1604.07149}{arXiv:1604.07149}.

\bibitem{TianyuMa}
Ma T., Geodesic rigidity of {L}evi-{C}ivita connections admitting essential
 projective vector fields,
 \href{https://doi.org/10.1007/s10711-019-00469-7}{\textit{Geom. Dedicata}}
 \textbf{205} (2020), 147--166,
 \href{https://arxiv.org/abs/1705.00460}{arXiv:1705.00460}.

\bibitem{MelnickFO}
Melnick K., Rigidity of transformation groups in differential geometry,
 \href{https://doi.org/10.1090/noti2279}{\textit{Notices Amer. Math. Soc.}}
 \textbf{68} (2021), 721--732,
 \href{https://arxiv.org/abs/2009.13937}{arXiv:2009.13937}.

\bibitem{MelnickNeusser2016}
Melnick K., Neusser K., Strongly essential flows on irreducible parabolic
 geometries, \href{https://doi.org/10.1090/tran/6814}{\textit{Trans. Amer.
 Math. Soc.}} \textbf{368} (2016), 8079--8110,
 \href{https://arxiv.org/abs/1410.4647}{arXiv:1410.4647}.

\bibitem{NaganoOchiai}
Nagano T., Ochiai T., On compact {R}iemannian manifolds admitting essential
 projective transformations, \textit{J.~Fac. Sci. Univ. Tokyo Sect.~IA Math.}
 \textbf{33} (1986), 233--246.

\bibitem{Sharpe1997}
Sharpe R.W., Differential geometry: {C}artan's generalization of {K}lein's
 {E}rlangen program, \textit{Grad. Texts in Math.}, Vol. 166, Springer, New
 York, 1997.

\end{thebibliography}
\end{document}